\documentclass[a4paper,10pt]{article}
\usepackage{amsmath,amssymb,amsfonts,amsthm}
\usepackage[all,cmtip]{xy}
\usepackage{arydshln}
\usepackage{a4wide}
\usepackage[usenames,dvipsnames]{color}
\usepackage[verbose,colorlinks=true,linktocpage=true,linkcolor=blue,citecolor=blue]{hyperref}
\usepackage{fouridx}
\usepackage{geometry}

\theoremstyle{definition}
\newtheorem{definition}{Definition}[section]

\theoremstyle{plain}
\newtheorem{theorem}[definition]{Theorem}

\newtheorem{proposition}[definition]{Proposition}
\newtheorem{lemma}[definition]{Lemma}

\theoremstyle{remark}
\newtheorem{remark}[definition]{Remark}
\theoremstyle{mystyle}
\newtheorem*{mainthm}{Main Theorem}

\numberwithin{equation}{section}

\allowdisplaybreaks[4]

\begin{document}
\title{Whittaker supermodules for the Takiff superalgebra $\mathfrak{osp}(1|2)$}
\author{Wei Wang ${}^1$, Yongjie Wang ${}^{2}$ and Chuchao Xu ${}^{3}$\footnote{Corresponding Author.}}
\maketitle

\begin{center}
\footnotesize
\begin{itemize}
\item[1] School of Mathematics and Statistics, Ningxia University, Yinchuan 750021, China.
\item[2] School of Mathematics, Hefei University of Technology, Hefei, Anhui, 230009, China.
\item[3] School of Mathematics, University of Science and Technology of China, Hefei, Anhui, 230026, China.
\end{itemize}
\end{center}

\begin{abstract} 
In this article, we classify all simple Whittaker supermodules over the Takiff Lie superalgebra \(\mathfrak g=\mathfrak{osp}(1|2)\otimes\mathbb C[t]/(t^2)\), associated with Whittaker characters defined on the even part of its positive nilpotent subalgebra.
\bigskip

\bigskip

\end{abstract}

\section{Introduction}

Whittaker modules constitute an important class of representations in the theory of Lie algebras, providing a natural framework for studying non-weight modules beyond the traditional highest-weight setting. Their study originated in the work of Arnal and Pinczon \cite{AP1974} on the Lie algebra $\mathfrak{sl}_2$ and was subsequently developed systematically by Kostant \cite{K1978} for finite-dimensional reductive Lie algebras. Since then, Whittaker modules and the corresponding categories have attracted considerable attention, particularly because of their connections with the Bernstein–Gelfand–Gelfand category $\mathcal O$ and finite $W$-algebras. These connections are exemplified by the Backelin functor \cite{Ba97} and the Skryabin equivalence \cite{Skr02}. Significant developments in these directions can be found, for example, in \cite{AB2021,Ba97,BM2011,MiSo97}. Further studies include modules induced from Whittaker modules \cite{M1985}, Whittaker modules for the Virasoro algebra \cite{OW2009}, simple Virasoro modules that are locally finite over a positive subalgebra \cite{MZ2013}, and singular Whittaker modules of type \(A_2\) \cite{XG2025}. These investigations demonstrate that Whittaker representations provide an effective approach to understanding the module categories of Lie algebras from a perspective complementary to highest-weight theory.

Whittaker modules and categories for various classes of Lie superalgebras have been investigated in \cite{BCW2014,C2021,CCh24b,CCh24a}, while related constructions have also been studied for infinite-dimensional Lie superalgebras; see, for instance, \cite{MGL2024,JXZh2024,LPX2019,WCM2026}. A fundamental distinction from the classical setting concerns the choice of Whittaker character. For a Lie superalgebra with a triangular decomposition, one may prescribe a character on the even part of its positive nilpotent subalgebra rather than imposing scalar relations on all positive root vectors. This distinction is particularly relevant when the square of an odd positive root vector is a nonzero even root vector. In such a situation, requiring the odd positive root vectors to annihilate a cyclic vector would impose additional restrictions on the possible eigenvalues of the even positive root vectors. Restricting the Whittaker character to the even part of the positive nilpotent subalgebra allows odd positive root vectors to act nontrivially on Whittaker vectors. Consequently, the structure of Whittaker supermodules depends not only on the Whittaker and central character, but also on the relations among the odd generators and the $\mathbb Z_2$-grading of the modules.

Another natural setting in which Whittaker theory exhibits new phenomena is provided by Takiff Lie algebras and their superanalogues. Originally introduced by Takiff \cite{T1971}, these algebras are obtained by extending a Lie algebra through a truncated polynomial ring. Their nonsemisimple structure produces additional abelian ideals and modifies the representation theory in ways that are not visible for the original Lie algebra \cite{CSZ2016, GM17, MS19}. Various aspects of the representation theory of Takiff superalgebras have been developed in recent years, including semisimple extensions, typical representations, and Whittaker modules; see \cite{CCS2026,CC2022}.

Recently, Chen and the second author of the article \cite{CW2025} investigated typical representations of Takiff Lie superalgebras and classified non-singular simple Whittaker modules in the Type I setting. However, the Kac-induction techniques used in these cases do not directly extend to Type II Takiff Lie superalgebras. Motivated by this, we investigate the Takiff Lie superalgebra \(\mathfrak g=\mathfrak{osp}(1|2)\otimes_{\mathbb C}\mathbb C[t]/(t^2)\), a basic example of Type II, and classify its simple Whittaker supermodules with respect to arbitrary characters of the even part of its positive nilpotent subalgebra.
The even part of the positive nilpotent subalgebra is
$$
\mathfrak n^+_{\bar 0}=\operatorname{Span}_{\mathbb C}\{e,\bar e\}.
$$

Given a Lie algebra homomorphism $\zeta:\mathfrak n^+_{\bar 0}\to\mathbb C$ with $\zeta(e)=a,\,\zeta(\bar e)=b$, we consider the $\mathfrak g$-supermodule generated by a nonzero Whittaker vector $w$ satisfying
$$
ew=aw,\qquad \bar ew=bw.
$$
The center of the universal enveloping algebra is the polynomial algebra $\mathbb C[z_1,z_2]$, with explicit central generators given in \cite{C2026}. Thus, the classification problem naturally involves both the Whittaker parameters $a,b$ and the central character $(\xi_1,\xi_2)$. Our main objective is to determine the structure of all simple Whittaker $\mathfrak g$-supermodules, describe their constructions, and establish their isomorphism classes.

We first investigate the case $a\ne 0$ and $b=0$. Here, the structure of simple Whittaker supermodules is essentially dependent on whether the first central element $z_1$ acts by zero. When $\xi_1=\xi_2=0$, $V$ is a non-singular simple Whittaker $\mathfrak{osp}(1|2)$-supermodule, which has already been classified by Chen and Cheng in \cite[Corollary 19]{CCh24a}. When $\xi_1=0$ and $\xi_2\ne 0$, we construct a simple Whittaker supermodule whose underlying ungraded module decomposes into a direct sum of two nonisomorphic simple modules (see Theorem \ref{3.8}). This provides an explicit illustration of how simplicity in the category of supermodules differs from simplicity in the category of ordinary modules. When $\xi_1\ne 0$, we construct a family of simple Whittaker supermodules parameterized by a square root $c$ of $\xi_1$. These constructions exhaust simple Whittaker supermodules with the prescribed central character and establish an isomorphism criterion that distinguishes the different choices of $c$ and parity. The resulting classification is stated in Theorems \ref{theo1} and \ref{3.17}.

Next, we consider the case $b\ne 0$. We first determine all Whittaker vectors in the universal Whittaker supermodule and construct an associative superalgebra $\mathcal A_\zeta$ that acts naturally on their space. This algebra is generated by two central elements and four odd elements satisfying Clifford-type relations. The universal Whittaker supermodule is shown to be a free right $\mathcal A_\zeta$-module, which allows us to recover simple Whittaker $\mathfrak g$-supermodules from simple $\mathcal A_\zeta$-supermodules by induction. More precisely, for every simple Whittaker supermodule $V$ of type $\zeta$, the space $\mathrm{Wh}_{\zeta}(V)$ of its Whittaker vectors is a simple $\mathcal A_\zeta$-supermodule, and
$$
V\cong M(\zeta)\otimes_{\mathcal A_{\zeta}} \mathrm{Wh}_{\zeta}(V)
$$
as $\mathfrak g$-supermodules (see Theorem \ref{Vphi2}).

Finally, we examine the zero Whittaker character $\zeta=0$. This case exhibits a close relationship between Whittaker representations, highest-weight representations, and modules over solvable Lie algebras. For $\zeta=0$, we show that every simple Whittaker supermodule is a highest-weight supermodule when $\xi_1\neq0$, while $\xi_1=0$ forces $\xi_2=0$. In the latter case, simple Whittaker supermodules either reduce to simple highest-weight supermodules over $\mathfrak{osp}(1|2)$ (see Theorem \ref{0osp}) or arise, up to parity reversal, by induction from simple modules over the 2-dimensional solvable Lie algebra $\mathfrak s=\mathbb Ch\oplus\mathbb C\bar f$ on which $\bar f$ acts injectively. Identifying $\mathrm U(\mathfrak s)$ with the Ore extension $\mathbb C[X][Y;X\frac{d}{dX}]$, we apply Bavula's classification \cite{Bavula2020} to describe these inducing modules and their isomorphism classes (see Theorem \ref{3.38}). The following is our main result.

\begin{mainthm}[Classification of simple Whittaker supermodules]
Let \(V\) be a simple Whittaker \(\mathfrak g\)-supermodule of type \(\zeta\). Then exactly one of the following three cases occurs:
\begin{itemize}
    \item[(1)] (Theorems \ref{3.8} and \ref{theo1}) Suppose that \(a\neq 0\) and \(b=0\). Then the following possibilities arise:
    \begin{enumerate}
        \item If \(\xi_1=\xi_2=0\), then \(V\) is a non-singular simple Whittaker \(\mathfrak{osp}(1|2)\)-supermodule, which has already been classified by Chen and Cheng in \cite[Corollary 19]{CCh24a}.
        \item If \(\xi_1=0\) and \(\xi_2\neq 0\), then
        \[
        V\cong \mathcal {S}\left(a,\frac{\xi_2}{2a}\right).
        \]
        \item If \(\xi_1\neq 0\), then there exist \(c\in\mathbb C^*\) with \(c^2=\xi_1\) and  \(\varepsilon\in\{0,1\}\), such that
       $$V\cong\Pi^\varepsilon \mathcal {S}(a,c,\xi_2).$$
    \end{enumerate}
    \item[(2)] (Theorem \ref{Vphi2}) Suppose that \(b\neq 0\). Then \(\mathrm{Wh}_{\zeta}(V)\) is a simple \(\mathcal A_{\zeta}\)-supermodule, and 
    $$V\cong M(\zeta)\otimes_{\mathcal A_{\zeta}}
    \mathrm{Wh}_{\zeta}(V).$$
    \item[(3)] (Theorems \ref{high} and \ref{3.38}) Suppose that \(a=b=0\). Then the following possibilities arise:
    \begin{enumerate}
        \item If \(\xi_1\neq 0\), then \(V\) is isomorphic to \(\Pi^{\varepsilon}\mathcal L(\lambda,\mu)\).
        \item If $\xi_1=0$ and $\xi_2\neq0$, then no simple Whittaker $\mathfrak g$-supermodule of type $\zeta=0$ exists.
        \item If \(\xi_1=\xi_2=0\), then $V$ is isomorphic to $\Pi^\varepsilon\left(\operatorname{Ind}_{\mathfrak a}^{\mathfrak g}S\right), \text{ for }\mathfrak pV\neq 0$, and to $ \Pi^{\varepsilon}\mathcal L(\lambda)$, otherwise.
        \end{enumerate}
\end{itemize}
The above cases exhaust all simple Whittaker \(\mathfrak g\)-supermodules. Their structures and isomorphism classes will be described explicitly in the subsequent classification theorems.
\end{mainthm}

Our results show that the classification of Whittaker supermodules over the Takiff extension of $\mathfrak{osp}(1|2)$ is governed by several distinct algebraic structures. When $b\neq 0$, the structure of the space of Whittaker vectors is described by a Clifford-type superalgebra. For certain singular characters, the distinction between ordinary modules and supermodules becomes essential, with parity reversal playing a significant role in the isomorphism classification. For the zero character, the problem reduces either to highest-weight theory or to the study of simple modules over a solvable Lie algebra. These classification results provide a systematic description of the simple Whittaker supermodules considered here and illustrate the additional representation-theoretic structures introduced by Takiff extensions and odd root operators.

The paper is organized as follows. In Section \ref{s2}, we recall the defining relations and triangular decomposition of the Takiff Lie superalgebra $\mathfrak g$, together with the description of the center of its universal enveloping algebra. Section \ref{s3} is devoted to the classification of simple Whittaker supermodules. In Subsections \ref{s31} and \ref{s32}, we treat the cases $a\ne0$, $b=0$ with $\xi_1=0$ and $\xi_1\ne0$, respectively. In Subsection \ref{s33}, we study the case $b\ne0$ through the associated Clifford-type superalgebra. Finally, Subsection \ref{s34} addresses the zero Whittaker character and its relationship with simple modules over Ore extensions.
\bigskip

\noindent{\bf Notations and terminologies:}
Throughout this paper, we will use the standard notations $\mathbb{Z}$, $\mathbb{Z}_+$, $\mathbb{N}$, $\mathbb{C}$ and $\mathbb{C}^*$ that represent the sets of integers, non-negative integers, positive integers, complex field and nonzero complex numbers, respectively. Denote by $\mathcal{Z}(\mathfrak{g})$ the center of the universal enveloping (super)algebra of a Lie (super)algebra $\mathfrak{g}$. The Kronecker delta $\delta_{ij}$ is equal to $1$ if $i=j$ and $0$ otherwise.

We write $\mathbb{Z}_2=\left\{\bar{0},\bar{1}\right\}$. For a homogeneous element $x$ of an associative or Lie superalgebra, we use $|x|$ to denote the degree of $x$ with respect to the $\mathbb{Z}_2$-grading. Throughout the paper, when we write $|x|$ for an element $x$,  we will always assume that $x$ is a homogeneous element and automatically extend the relevant formulas by linearity (whenever applicable). All modules of Lie superalgebras are assumed to be $\mathbb{Z}_2$-graded. Unless otherwise specified, $\otimes$ always denotes the tensor product over $\mathbb{C}$.

\medskip
\noindent
{\bf Acknowledgement}. The paper is partially supported by the National Natural Science Foundation of China (Nos. 12361006, 12471025), and Y. Wang is supported by the Anhui Provincial Natural Science Foundation 2308085MA01, and the Fundamental Research Funds for the Central Universities JZ2025HGTG0251. The authors thank Yan-an Cai, Chih-Whi Chen, and Shun-Jen Cheng for interesting discussions and helpful comments.

\section{Takiff superalgebra}\label{s2}
The complex Lie superalgebra $\mathfrak{osp}(1|2)=\mathfrak{osp}(1|2)_{\bar{0}}\oplus\mathfrak{osp}(1|2)_{\bar{1}}$ admits a basis $\{e,f,h,E,F\}$, where
$$\mathfrak{osp}(1|2)_{\bar{0}}={\rm Span}_{\mathbb{C}}\{e,f,h\},\ \ 
\mathfrak{osp}(1|2)_{\bar{1}}={\rm Span}_{\mathbb{C}}\{E,F\},$$
whose non-trivial superbrackets are given by
\begin{eqnarray}
\label{osp1}&&[e,f]=h,\ \ [h,e]=2e,\ \ [h,f]=-2f,\\
\label{osp2}&&[E,E]=2e,\ \ [F,F]=-2f,\ \ [E,F]=h\\
\label{osp3}&&[e,F]=-E,\ \ [f,E]=-F,\ \ [h,E]=E,\ \ [h,F]=-F.
\end{eqnarray}

What we consider in the paper is the non-semisimple Lie superalgebra $\mathfrak{g}=\mathfrak{osp}(1|2)\otimes\mathbb{C}[t]/(t^2)$ with the brackets
 $$[a\otimes t^i,b\otimes t^j]=[a,b]\otimes t^{i+j},$$
 where $a,b\in\{e,f,h,E,F\}$ and $i,j\in\{0,1\}$, with $t$ even, $t^2=0$, and the bracket on the right hand side determined by \eqref{osp1}-\eqref{osp3}. Thus, $\mathfrak{g}$ is the first Takiff extension of $\mathfrak{osp}(1|2)$, called the {\it $1$-Takiff superalgebra of type $\mathfrak{osp}(1|2)$}. We identify $x\in\mathfrak{osp}(1|2)$ with $x\otimes 1$ and write $\bar x=x\otimes t$.
 
The Takiff Lie superalgebra has superdimension $(6|4)$, with 
 $$\mathfrak{g}_{\bar{0}}={\rm Span}_{\mathbb{C}}\{e,f,h,\bar{e},\bar{f},\bar{h}\},\quad \mathfrak{g}_{\bar{1}}={\rm Span}_{\mathbb{C}}\{E,F,\bar{E},\bar{F}\},$$
admitting the following relations
$$[\bar{x},\bar{y}]=0,\quad [x,\bar{y}]=[\bar{x},y]=\overline{[x,y]},\quad x,y\in\{e,f,h,E,F\}.$$
The triangular decomposition of $\mathfrak{g}$ is 
$$\mathfrak{g}=\mathfrak{n}^{-}\oplus\mathfrak{h}\oplus\mathfrak{n}^{+},$$
with 
$$\mathfrak{n}^{-}={\rm Span}_{\mathbb{C}}\{f,\bar{f},F,\bar{F}\} ,\quad\mathfrak{h}={\rm Span}_{\mathbb{C}}\{h,\bar{h}\},\quad\mathfrak{n}^{+}={\rm Span}_{\mathbb{C}}\{e,\bar{e},E,\bar{E}\}.$$

For homogeneous elements, the bracket in the enveloping algebra is the supercommutator $[x,y]=xy-(-1)^{|x||y|}yx$. In particular, $E^2=e,F^2=-f$ and $\bar E^2=\bar F^2=0$.
The following identities can be established by induction.
\begin{lemma}\label{[]}
    Let $n\geq 0$, $\bar x\in\{\bar e,\bar f,\bar h,\bar E,\bar F\}$ and $y\in\{e,f,h,E,F,\bar e,\bar f,\bar h,\bar E,\bar F\}.$ Then the supercommutators $[\bar x,y^n]$ are given by the following formulas:
\begin{itemize}
    \item[(1).] For \(y\in \{e,f\}\), we have
\begin{gather*}
[\bar e,e^n]=[\bar f,f^n]=[\bar E,e^n]=[\bar F,f^n]=0,\\ 
[\bar f,e^n]=-ne^{n-1}\bar h-n(n-1)e^{n-2}\bar e,\qquad [\bar e,f^n]=nf^{n-1}\bar h-n(n-1)f^{n-2}\bar f,\\
[\bar h,e^n]=2ne^{n-1}\bar e,\qquad [\bar h,f^n]=-2nf^{n-1}\bar f,\qquad [\bar E,f^n]=nf^{n-1}\bar F,\qquad [\bar F,e^n]=ne^{n-1}\bar E.
\end{gather*}
\item[(2).] For \(y=h\), we have
\[
[\bar x,h^n]
=\bigl((h-\lambda_{\bar x})^n-h^n\bigr)\bar x,
\]

where
\[
\lambda_{\bar e}=2,\quad
\lambda_{\bar f}=-2,\quad
\lambda_{\bar h}=0,\quad
\lambda_{\bar E}=1,\quad
\lambda_{\bar F}=-1.
\]
\item[(3).] For \(y\in \{E,F\}\), write \(n=2m\) or \(n=2m+1\) for any $m\in\mathbb{Z_+}$. Then
\[
[\bar x,E^{2m}]=[\bar x,e^m],
\qquad
[\bar x,F^{2m}]=(-1)^m[\bar x,f^m].
\]

For \(m\in\mathbb{Z}_+\), the odd powers satisfy
\begin{alignat*}{2}
[\bar{e}, E^{2m+1}] =& 0,\qquad &&[\bar{f}, F^{2m+1}] = 0, \\[1ex]
[\bar{h}, E^{2m+1}] =& 2m e^{m-1} \bar{e} E + e^m \bar{E},\qquad && [\bar{h}, F^{2m+1}] = (-1)^m (-2m f^{m-1} \bar{f} F - f^m \bar{F}), \\[1ex]
[\bar{E}, E^{2m+1}] =& 2e^m \bar{e},\qquad && [\bar{F}, F^{2m+1}] = 2(-1)^{m+1} f^m \bar{f}, \\[1ex]
[\bar{F}, E^{2m+1}] =& m e^{m-1} \bar{E} E + e^m \bar{h},\qquad &&[\bar{E}, F^{2m+1}] = (-1)^m (m f^{m-1} \bar{F} F + f^m \bar{h}),
\end{alignat*}
together with 
\begin{align*}
[\bar{f}, E^{2m+1}] =& -m e^{m-1} \bar{h} E - m(m-1) e^{m-2} \bar{e} E - e^m \bar{F},\\[1ex]
[\bar{e}, F^{2m+1}] =& (-1)^m (m f^{m-1} \bar{h} F - m(m-1) f^{m-2} \bar{f} F - f^m \bar{E}).
\end{align*}
Terms involving negative powers are understood to vanish.
\item[(4).] For \(y\in\{\bar e,\bar f,\bar h,\bar E,\bar F\}\), we have
$[\bar x,y^n]=0.$
\end{itemize}
\end{lemma}
\begin{definition} 
Let $\mathfrak n_{\bar 0}^+={\rm Span}_{\mathbb C}\{e,\bar e\}$ and $\zeta:\mathfrak n_{\bar 0}^{+}\to \mathbb C$ be a Lie algebra homomorphism. Set $a=\zeta(e),b=\zeta(\bar e)$. For a $\mathfrak g$-supermodule $V$, we define 
\begin{itemize}
    \item $\mathrm{Wh}_{\zeta}(V)=\{v\in V|(e-a)v=(\bar e-b)v=0\}$. A nonzero element of  $\mathrm{Wh}_{\zeta}(V)$ is called a \textit{Whittaker vector} of type $\zeta$.
    \item $V$ is called a \textit{Whittaker supermodule} of type $\zeta$ if it is generated by a Whittaker vector.
    \item Define the \textit{Universal Whittaker supermodule} of type $\zeta$ by $M(\zeta):=\mathrm U(\mathfrak g)\otimes_{\mathrm U(\mathfrak n_{\bar 0}^+)}\mathbb Cv_{\zeta}$, where $v_{\zeta}$ is even with $ev_{\zeta}=av_{\zeta}$ and $\bar ev_{\zeta}=bv_{\zeta}$.
\end{itemize}
\end{definition}
For convenience, we define the following subalgebras of $\mathfrak{g}$:
\begin{align*}
     \mathfrak{a}=&{\rm Span}_\mathbb{C}\{h, \bar{h}, e,\bar{f}, \bar{e},E,\bar{E}, \bar{F}\},
    &\mathfrak{k}=&{\rm Span}_\mathbb{C}\{\bar{h}, e,\bar{e},E, \bar{E}\},\qquad
   \mathfrak{l}={\rm Span}_\mathbb{C}\{\bar{h},e,\bar{e},\bar{E},\bar{F}\},\\
    \mathfrak{p}=&{\rm Span}_\mathbb{C}\{\bar{h},\bar{e},\bar{f},\bar{E},\bar{F}\}, &\mathfrak{s}=&{\rm Span}_\mathbb{C}\{h,\bar{f}\}.
\end{align*}
\begin{remark}
Let $\psi: \mathfrak p\to  \mathbb{C}$ be a Lie superalgebra homomorphism. A nonzero vector $v\in V$ is called a \textit{quasi-Whittaker vector} of type $\psi$ if $xv=\psi(x)v$ for all $x\in \mathfrak p$. The $\mathfrak g$-module $V$ is called a \textit{quasi-Whittaker module} of type $\psi$ if it is generated by such a vector. Simple quasi-Whittaker modules for the Takiff Lie superalgebra of $\mathfrak {osp}(1|2)$ were classified in \cite{C2026}.
\end{remark}

\begin{proposition}[{\cite[Theorem 3.1]{C2026}}]{\label{center}}
The center of $\mathrm{U}(\mathfrak{g})$ is $\mathcal{Z}(\mathfrak{g})=\mathbb{C}[z_1,z_2]$, where
\begin{eqnarray*}
z_1=4\bar{f}\bar{e}+\bar{h}^2+2\bar{F}\bar{E},\ \ z_2=2{f}\bar{e}+h\bar{h}+2{e}\bar{f}+{F}\bar{E}-E\bar{F}.
\end{eqnarray*}
\end{proposition}
\section{Whittaker supermodules}\label{s3}
Throughout this section, let $V$ be a simple Whittaker $\mathfrak g$-supermodule of type $\zeta$, generated by a Whittaker vector $w$. Since $\mathrm U(\mathfrak g)$ is countable-dimensional over $\mathbb C$, the super version of Schur's lemma implies that $z_1,z_2$ act on $V$ by scalars $\xi_1,\xi_2\in \mathbb C$, respectively. Thus,
\begin{equation}\label{ca}
    z_1w=\xi_1w,\qquad z_2w=\xi_2w.
\end{equation}
\begin{remark}
    If $V$ is generated by a Whittaker vector $w$, then each nonzero homogeneous component of $w$ remains a Whittaker vector (since $e$ and $\bar e$ are even elements). Therefore, we can choose a homogeneous Whittaker generator.
\end{remark}
\begin{lemma}\label{cons}
    Let $V$ be a simple Whittaker $\mathfrak g$-supermodule with a homogeneous Whittaker vector $w$ of type $\zeta$. Suppose that there exists a homogeneous element $x\in \mathrm U(\mathfrak g)$ such that $x^2w=0$ and $xw\in \mathrm{Wh}_{\zeta}(V)$. Then $V$ contains a nonzero homogeneous cyclic Whittaker vector 
    \begin{equation}\label{nv}
        v=
    \begin{cases}
    w, &\text{ if } xw=0,\\
    xw,    &\text{ if } xw\neq0,

    \end{cases}
    \end{equation}
such that $xv=0$.
\end{lemma}
\begin{proof}
    If $xw=0$, take $v=w$. Otherwise, take $v=xw$. In the latter case, $v$ is nonzero and homogeneous, belongs to $\mathrm{Wh}_{\zeta}(V)$ by assumption, and satisfies $xv=x^2w=0$. Since $V$ is simple, the nonzero subsupermodule $\mathrm U(\mathfrak g)v$ equals $V$. Thus, \(v\) is a homogeneous cyclic Whittaker vector satisfying $xv=0$.
\end{proof}

In the following subsections, we distinguish the following four cases: 
$$(1).\ a\neq0,\,b=0,\,\xi_1=0;\quad (2).\ a\neq 0,\,b=0,\,\xi_1\neq 0;\quad (3).\ b\neq 0;\quad (4).\ a=b=0.$$
\subsection{Singular Whittaker supermodules with $a\neq0,b=0,\xi_1=0$}\label{s31}
Assume $a\neq 0,b=0$ and $\xi_1=0$. We first construct a homogeneous cyclic Whittaker vector satisfying additional annihilation relations.
\begin{proposition}\label{xi0}
    There exists a nonzero homogeneous cyclic Whittaker vector $\widetilde w\in V$ such that
    \begin{equation}\label{aneq0}
        e\widetilde w=a\widetilde w,\qquad \bar e\widetilde w=\bar E\widetilde w=\bar h\widetilde w=\bar F\widetilde w=0,\qquad \bar f\widetilde w=\lambda \widetilde w,
    \end{equation}
    where $\lambda=\frac{\xi_2}{2a}$.
\end{proposition}
\begin{proof}
    Let $w_0$ be a nonzero homogeneous cyclic Whittaker vector of $V$. Since $\bar E^2w_0=0$ and $\bar Ew_0\in \mathrm{Wh}_{\zeta}(V)$, by Lemma \ref{cons}, we have a nonzero homogeneous cyclic Whittaker vector $v_1$ satisfying $\bar Ev_1=0$. Furthermore, by \eqref{ca} and Proposition \ref{center}, we have $z_1v_1=\bar h^2v_1=0$. The defining relations give $\bar hv_1\in \mathrm{Wh}_{\zeta}(V)$. The Lemma \ref{cons} also produces a nonzero homogeneous cyclic Whittaker vector $v_2$ with $\bar hv_2=0$. The construction in \eqref{nv} also gives $\bar Ev_2=0$. 

    Obviously, we have $(\bar FE)^2v_2=\bar F\bar hEv_2=\bar F(E\bar h+\bar E)v_2=0$ and $(\bar FE)v_2\in \mathrm{Wh}_{\zeta}(V)$. By Lemma \ref{cons}, there exists a nonzero homogeneous cyclic Whittaker vector $v_3$ satisfying $(\bar FE)v_3=0$. By the construction in \eqref{nv}, we have $\bar Ev_3=\bar hv_3=0$. Moreover, we have 
    $$
    0=E(\bar FE)v_3=(\bar h-\bar FE)Ev_3=(-a\bar F+E\bar h+\bar E)v_3=-a\bar Fv_3.
    $$
    Since $a\neq 0$, we have $\bar Fv_3=0$. Furthermore, by \eqref{ca}, we obtain
    $$
    \xi_2v_3=z_2v_3=(2a\bar f+\bar FE)v_3=2a\bar fv_3.
    $$
    Taking $\widetilde w=v_3$ completes the proof.
\end{proof}
\begin{remark}\label{toeven}
    The vector \(\widetilde w\) in Proposition \ref{xi0} can be chosen even. If $\widetilde w$ is odd, then $E\widetilde w$ is even and nonzero because $E^2\widetilde w=a\widetilde w\ne0$. Moreover, \(E\widetilde w\) remains cyclic and satisfies \eqref{aneq0}. Replacing \(\widetilde w\) by \(E\widetilde w\) therefore gives the required even vector.
\end{remark}

If \(\xi_2=0\), then \(\lambda=0\), and \eqref{aneq0} implies $\mathfrak p\widetilde w=0$. Since \(\mathfrak p\) is an ideal of \(\mathfrak g\) and \(V=\mathrm U(\mathfrak g)\widetilde w\), it follows that \(\mathfrak pV=0\). Hence, \(V\) factors through \(\mathfrak{osp}(1|2)\). Moreover, the supermodules
\[
\{L(\mu,\zeta)\mid \mu(h)\notin\mathbb Z_{\leqslant 0}\}
\]
are all simple Whittaker \(\mathfrak{osp}(1|2)\)-supermodules; see \cite[Corollary 19]{CCh24a} for notation $L(\mu,\zeta)$ and other more details.

We therefore assume $\xi_2\neq 0$ throughout the remainder of this subsection.

\begin{definition}
    Let $w_{\zeta}$ be an even free cyclic generator in the following presentation. Define the quotient $\mathfrak g$-supermodule  by
    $$
    \mathcal {S}(a,\lambda)=\mathrm U(\mathfrak g)w_{\zeta}/\langle ew_{\zeta}-aw_{\zeta},\bar ew_{\zeta},\bar hw_{\zeta},\bar Ew_{\zeta},\bar Fw_{\zeta},\bar fw_{\zeta}-\lambda w_{\zeta}\rangle
    $$ 
    and denote the image of $w_{\zeta}$ by $\bar w_{\zeta}$. Choose $s\in \mathbb C$ such that $s^2=a$. For $\varepsilon\in \{-1,1\}$, define the quotient $\mathfrak g$-module 
    $$
    \mathcal {S}_{\varepsilon s}(a,\lambda)=\mathrm U(\mathfrak g)w_{\zeta}/\langle ew_{\zeta}-aw_{\zeta},Ew_{\zeta}-\varepsilon sw_{\zeta},\bar ew_{\zeta},\bar hw_{\zeta},\bar Ew_{\zeta},\bar Fw_{\zeta},\bar fw_{\zeta}-\lambda w_{\zeta}\rangle
    $$
    with cyclic vectors $v_{\varepsilon s}$. 
\end{definition}
    Then PBW theorem shows that the following sets form bases of $\mathcal {S}(a,\lambda)$ and $\mathcal {S}_{\varepsilon s}(a,\lambda)$, respectively.
    $$
    \{f^mh^nF^{\alpha}E^{\beta}\bar w_{\zeta}\mid m,n\geq 0,\alpha,\beta  \in \{0,1\}\},\qquad \{f^mh^nF^{\alpha}v_{\varepsilon s}|m,n\geq 0,\alpha  \in \{0,1\}\}.
    $$ 

    We first establish the simplicity of $\mathcal {S}_{\varepsilon s}(a,\lambda)$.
\begin{lemma}\label{xisi}
    If $a\ne 0$ and $\lambda\neq0$, then $\mathcal {S}_{\varepsilon s}(a,\lambda)$ is a simple $\mathfrak g$-module for $\varepsilon=\pm 1$. 
\end{lemma}
\begin{proof}
    Let $N$ be a nonzero submodule of $\mathcal {S}_{\varepsilon s}(a,\lambda)$. It suffices to show that \(v_{\varepsilon s}\in N\).
    
    Let $0\neq v\in N$. We write $v=P_0(f,h)v_{\varepsilon s}+P_1(f,h)Fv_{\varepsilon s}$, where $P_0(f,h),P_1(f,h)\in \mathbb C[f,h]$. Then $\bar Fv=-2\lambda P_1(f,h+1)v_{\varepsilon s}\in N$. 
    If \(P_1(f,h)=0\), then \(P_0(f,h)v_{\varepsilon s}\in N\). Otherwise, \(\lambda\ne0\) implies that \(P_1(f,h+1)v_{\varepsilon s}\in N\). In either case, \(N\) contains a nonzero vector \(v_1=P(f,h)v_{\varepsilon s}\) with \(P(f,h)\in\mathbb C[f,h]\).
    
    Write $P(f,h)=\sum_{i=0}^{d_0}f^iq_i(h)$, where $q_{d_0}(h)\neq 0$. Then, $\bar h^{d_0}v_1=(-2\lambda)^{d_0}d_0!q_{d_0}(h+2d_0)v_{\varepsilon s}\in N$. Since $\lambda\neq 0$, we have $q_{d_0}(h+2d_0)v_{\varepsilon s}\in N$. 
    
    Thus, assume $v_2=q(h)v_{\varepsilon s}\in N$, where $q(h)\in \mathbb C[h]$ is nonzero. Let $n=\deg q(h)$ and $q_n\neq 0$ be its leading coefficient. Then, $(\bar f-\lambda)^nv_2=(2\lambda)^nn!q_nv_{\varepsilon s}\in N$. As $\lambda\neq 0$, this implies $v_{\varepsilon s}\in N$. Hence, $\mathcal {S}_{\varepsilon s}(a,\lambda)$ is a simple $\mathfrak g$-module.
\end{proof}
\begin{remark}\label{ves}
    The proof of Lemma \ref{xisi} also shows that 
    $$\{v\in \mathcal {S}_{\varepsilon s}(a,\lambda)|\bar Fv=\bar hv=(\bar f-\lambda)v=0\}=\mathbb Cv_{\varepsilon s}.$$
\end{remark}
\begin{lemma}\label{ncong}
    If $a\ne0$ and $\lambda\neq 0$, then $\mathcal {S}_{s}(a,\lambda)$ and $\mathcal {S}_{-s}(a,\lambda)$ are nonisomorphic $\mathfrak g$-modules.
\end{lemma}
\begin{proof}
    Suppose, to the contrary, that there exists a $\mathfrak g$-isomorphism
    $$
    \sigma:\mathcal {S}_{s}(a,\lambda)\to \mathcal {S}_{-s}(a,\lambda).
    $$
    As $\sigma(v_{s})\neq 0$ and $\bar Fv_{s}=\bar hv_{s}=(\bar f-\lambda)v_{s}=0$, we obtain
    \begin{align*}
        \bar F\sigma(v_s)=\bar h\sigma(v_s)=(\bar f-\lambda)\sigma(v_s)=0.
    \end{align*}
    Thus, by Remark \ref{ves}, there exists a scalar $c\in \mathbb C^*$ such that $\sigma(v_{s})=cv_{-s}$. On the other hand, we have 
    $$
    scv_{-s}=s\sigma(v_{s})=\sigma(sv_s)=\sigma(Ev_s)=E\sigma(v_{s})=cEv_{-s}=-scv_{-s},
    $$
    which is a contradiction, since $s,c\neq 0$. Therefore, $\mathcal {S}_{s}(a,\lambda)$ and $\mathcal {S}_{-s}(a,\lambda)$ are nonisomorphic $\mathfrak g$-modules.
\end{proof}
We now use Lemmas \ref{xisi} and \ref{ncong} to classify simple Whittaker $\mathfrak g$-supermodules in this case.
\begin{theorem}\label{3.8}
Suppose that $a\ne 0,\xi_2\neq0,\xi_1=0$ and $s^2=a$. Let $\lambda=\frac{\xi_2}{2a}$. Then the following statements hold:
\begin{enumerate}
    \item[(1).] As a $\mathfrak g$-module, $\mathcal {S}(a,\lambda)\cong \mathcal {S}_{s}(a,\lambda)\oplus \mathcal {S}_{-s}(a,\lambda)$.
    \item[(2).] \(\mathcal {S}(a,\lambda)\) is a simple $\mathfrak g$-supermodule.
    \item[(3).] Every simple Whittaker $\mathfrak g$-supermodule \(V\) of type \(\zeta\) with central character \((\xi_1,\xi_2)\) is isomorphic to \(\mathcal {S}(a,\lambda)\) as a $\mathfrak g$-supermodule.
    \item[(4).] $\mathcal {S}(a,\lambda)\cong \Pi \mathcal {S}(a,\lambda)$.
\end{enumerate}
\end{theorem}
\begin{proof}
   (1) Define $w_+=(E+s)\bar w_{\zeta},w_-=(E-s)\bar w_{\zeta}$. The PBW basis ensures that \(w_+\) and \(w_-\) are nonzero. A direct calculation gives
    \begin{equation}\label{pm}
        \begin{gathered}
        ew_\pm=aw_\pm,\quad
        Ew_\pm=\pm sw_\pm,\quad  \bar fw_\pm=\lambda w_\pm,\\
        \bar ew_\pm=\bar hw_\pm
        =\bar Ew_\pm=\bar Fw_\pm=0.
        \end{gathered}
    \end{equation}
    Let $W_+=\mathrm U(\mathfrak g)w_+$ and $W_-=\mathrm U(\mathfrak g)w_-$. The relations induce surjective $\mathfrak g$-module homomorphisms
    \begin{alignat*}{2}
        \mathcal {S}_{s}(a,\lambda) &\twoheadrightarrow W_+, \qquad& v_s &\mapsto w_+, \\
        \mathcal {S}_{-s}(a,\lambda) &\twoheadrightarrow W_-, \qquad& v_{-s} &\mapsto w_-.
    \end{alignat*}
    By Lemma \ref{xisi}, $\mathcal {S}_{s}(a,\lambda)$ and $\mathcal {S}_{-s}(a,\lambda)$ are simple $\mathfrak g$-modules. Since both $W_+$ and $W_-$ are nonzero, the two homomorphisms are isomorphisms. Hence
    \begin{equation}\label{lw}
    \mathcal {S}_{s}(a,\lambda)\cong W_+,
    \qquad
    \mathcal {S}_{-s}(a,\lambda)\cong  W_-.
    \end{equation}
    Moreover, we have $\bar w_{\zeta }=\frac{w_+-w_-}{2s}$ and $E\bar w_{\zeta}=\frac{w_++w_-}{2}$. It follows that $\bar w_{\zeta}\in W_++W_-$. Since \(\bar w_{\zeta}\) generates \(\mathcal {S}(a,\lambda)\), we have $\mathcal {S}(a,\lambda)=W_++W_-$. If \(W_+\cap W_-\neq0\), then the simplicity of $W_+$ and $W_-$ implies that $W_+\cap W_-=W_+=W_-$. Consequently, $\mathcal {S}_{s}(a,\lambda)\cong \mathcal {S}_{-s}(a,\lambda)$. This contradicts Lemma \ref{ncong}. Therefore,
    \begin{equation}\label{op}
        \mathcal {S}(a,\lambda)=W_+\oplus  W_-\cong \mathcal {S}_{s}(a,\lambda)\oplus  \mathcal {S}_{-s}(a,\lambda)
    \end{equation}
    as $\mathfrak g$-modules.
    
    (2) Let \(P\) denote the parity operator on \(\mathcal {S}(a,\lambda)\), defined on homogeneous vectors by \(P(v)=(-1)^{|v|}v\). Then
    \begin{equation}\label{change}
        P(W_+)=W_-,\qquad P(W_-)=W_+.       
    \end{equation}

    Let $N$ be a nonzero $\mathfrak g$-subsupermodule of $\mathcal {S}(a,\lambda)$. By \eqref{op} and Lemma \ref{ncong}, $\mathcal {S}(a,\lambda)$, viewed as a $\mathfrak g$-module, is a direct sum of two nonisomorphic simple $\mathfrak g$-modules. Therefore, its $\mathfrak g$-submodules are precisely 
    $$
    0,\qquad W_+,\qquad W_-,\qquad  \mathcal {S}(a,\lambda).
    $$
    
    On the other hand, since $N$ is a $\mathfrak g$-supermodule, it satisfies $P(N)=N$. However, by \eqref{change}, $P(W_+)=W_-\neq W_+$ and $P(W_-)=W_+ \neq W_-$. Thus, neither \(W_+\) nor \(W_-\) is a $\mathfrak g$-subsupermodule. Since \(N\neq0\), we have $N=\mathcal {S}(a,\lambda)$. Therefore, $\mathcal {S}(a,\lambda)$ is a simple $\mathfrak g$-supermodule.

    (3) By Proposition \ref{xi0} and Remark \ref{toeven}, $V$ contains a nonzero even cyclic Whittaker vector $\widetilde w$ satisfying \eqref{aneq0}.

    Define a map
    \[
        \pi:\mathcal {S}(a,\lambda)\to V,\qquad x\bar w_{\zeta}\mapsto x\widetilde w,\quad \text{ for } x\in \mathrm U(\mathfrak g).
    \]
    The defining relations show that \(\pi\) is a well-defined even $\mathfrak g$-supermodule homomorphism. Since \(\widetilde w\) is cyclic, \(\pi\) is surjective. Its kernel is zero by the simplicity of \(\mathcal {S}(a,\lambda)\), and therefore \(V\cong \mathcal {S}(a,\lambda)\) as $\mathfrak g$-supermodules.

    (4) Since $\Pi \mathcal {S}(a,\lambda)$ is also a simple Whittaker $\mathfrak g$-supermodule of type $\zeta$ with central character $(\xi_1,\xi_2)$, it follows from (3) that \(\mathcal {S}(a,\lambda)\cong\Pi \mathcal {S}(a,\lambda)\).
\end{proof}
\subsection{Singular Whittaker supermodules with $a\neq 0,b=0$, and $\xi_1\neq 0$}\label{s32}
Assume $a\neq 0,b=0$, and $\xi_1\neq 0$. Choose $s,c\in \mathbb C^*$ with $s^2=a$ and $c^2=\xi_1$. Define the 1-dimensional (ungraded) $\mathfrak k$-module $\mathbb Cw_{\zeta}$ by 
$$
ew_{\zeta}=aw_{\zeta},\ Ew_{\zeta}=sw_{\zeta},\ \bar hw_{\zeta}=cw_{\zeta},\ \bar ew_{\zeta}=\bar Ew_{\zeta}=0.
$$
Define the induced $\mathfrak g$-module and its central quotient
$$
W=\mathrm U(\mathfrak g)\otimes_{\mathrm U(\mathfrak k)}\mathbb Cw_{\zeta},\qquad L(s,c,\xi_2)=W/(z_2-\xi_2)W.
$$
We regard both $W$ and $L(s,c,\xi_2)$ as ungraded $\mathfrak g$-modules.
Denote the image of $w_{\zeta}$ in the quotient by $\bar w_{\zeta}$. Proposition \ref{center} and $z_2\bar w_{\zeta}=\xi_2\bar w_{\zeta}$ give
\begin{equation*}
    \bar f\bar w_{\zeta}=\frac{\xi_2-c(h+1)}{2a}\bar w_{\zeta}-\frac{\bar F}{2s}\bar w_{\zeta}.
\end{equation*}

We first show that $L(s,c,\xi_2)$ is nonzero and simple.
\begin{lemma}\label{Lbas}
    For any \(s,c\in\mathbb C^*\) and \(\xi_2\in\mathbb C\), the $\mathfrak g$-module \(L(s,c,\xi_2)\) is nonzero. Moreover, it admits a basis consisting of the vectors
\[
\mathcal B=\left\{
f^mh^nF^\varepsilon\bar F^\delta \bar w_{\zeta}
\;\middle|\;
m,n\ge0,
\varepsilon,\delta\in\{0,1\}
\right\}.
\]
\end{lemma}
\begin{proof}
    By the PBW theorem, $W$ has a basis $\{f^m h^n F^\varepsilon
    \bar f^k \bar F^\delta w_{\zeta}|m,n,k\geq 0,\varepsilon,\delta \in \{0,1\}\}$. In particular, $W\neq 0$. Let \(\deg_{\bar f}\) denote the \(\bar f\)-degree relative to this basis. Since \(z_2\) is central, we obtain
    \begin{eqnarray}\label{tri}
    z_2(f^m h^n F^\varepsilon\bar f^k\bar F^\delta w_{\zeta})=
    2a f^m h^n F^\varepsilon \bar f^{k+1}\bar F^\delta w_{\zeta}
    +\text{terms of degree at most }k
    \text{ in }\bar f.
    \end{eqnarray}
    Since \(a\neq0\), we have 
    \[
    \deg_{\bar f}\bigl((z_2-\xi_2)v\bigr)
    =
    \deg_{\bar f}v+1,
    \]
    for any nonzero $v\in W$. Consequently, no nonzero vector of \(\bar f\)-degree zero belongs to \((z_2-\xi_2)W\). Thus, the images of the $\bar f$-degree zero PBW basis vectors are linearly independent in \(L(s,c,\xi_2)\).

    To prove spanning, let $v\in W$ have $\deg_{\bar f}(v)=k\geq 1$. Equation \eqref{tri} allows us to eliminate its leading $\bar f$-degree terms subtracting \((z_2-\xi_2)u\) for a suitable \(u\in W\) with $\deg_{\bar f}u=k-1$. Repeating this reduction yields a representative of \(v+(z_2-\xi_2)W\) of $\bar f$-degree zero. 
    Hence, \(\mathcal B\) is a basis of \(L(s,c,\xi_2)\). In particular, $L(s,c,\xi_2)\neq0$.

\end{proof}
\begin{lemma}\label{sim2}
    For any $s,c\in \mathbb C^*$ and $\xi_2\in \mathbb C$, the $\mathfrak g$-module $L(s,c,\xi_2)$ is a simple $\mathfrak g$-module.
\end{lemma}
\begin{proof}
     By Lemma \ref{Lbas}, the vectors $f^m h^n F^\varepsilon
     \bar F^\delta \bar w_{\zeta}$, with $m,n\geq 0$ and $\varepsilon,\delta\in\{0,1\}$, form a basis of $L(s,c,\xi_2)$. Set
    $$
    x=\bar F\bar w_{\zeta},\ \ H=\mathbb C[h]\bar w_{\zeta}\oplus \mathbb C[h]x,\ \ K=H\oplus FH.
    $$
    Then we have the decomposition 
    \begin{equation}\label{struL}
        L(s,c,\xi_2)=\bigoplus\limits_{m\geq 0}f^mK.
    \end{equation}
    Let $N$ be a nonzero submodule of $L(s,c,\xi_2)$, and choose a nonzero vector $v\in N$. By \eqref{struL}, we write $v=\sum\limits_{m=0}^{d_1}f^mk_m$, where $k_m\in K$ and $k_{d_1}\neq 0$.
    Note that
    \begin{equation}\label{actbare}
        \bar ef^mk_m=mf^{m-1}\bar h k_m-m(m-1)f^{m-2}\bar fk_m.
    \end{equation}
    On $K$, the operator $\bar h$ has the form $\bar h|_{K}=c\text{id}_{K}+T$, where $T$ is a linear operator on $K$ satisfying $T^2=0$. Since $c\neq 0$, the operator $\bar h$ is invertible on $K$. By \eqref{actbare} and the fact that $\bar hK,\bar fK\subseteq K$, we have  $ \bar e^{d_1}v=d_1!\bar h^{d_1}k_{d_1}\neq 0$. Hence, $N\cap K\neq 0$. Let $0\neq v_1\in N\cap K$ and write 
    $v_1=h_0+Fh_1$ for some $h_0,h_1\in H$. Then
    $$
        \bar Ev_1=\bar Eh_0+\bar EFh_1=\bar hh_1=ch_1.
    $$
    If $h_1=0$, then $v_1\in N\cap H$ is nonzero. Otherwise, $\bar Ev_1=ch_1\in N\cap H$ is nonzero because $c\neq 0$. Thus, in either case, $N\cap H\neq 0$. Choose a nonzero vector $0\neq v_2=P(h)\bar w_{\zeta}+Q(h)x\in N\cap H$, where $P(h),Q(h)\in \mathbb C[h]$. We have 
    \begin{equation}\label{abarF}
        \bar Fv_2=\bar F(P(h)\bar w_{\zeta }+Q(h)x)=P(h+1)x.
    \end{equation}

    Suppose that $\text{deg}_hP(h)=n$ and let $p_n\neq 0$ denote its leading coefficient‌. Then 
    \begin{equation}\label{ae-a}
        (e-a)^nP(h)x=(-2a)^nn!p_nx\neq 0.
    \end{equation}
    If $P(h)\neq 0$, then \eqref{abarF} gives $0\neq P(h+1)x\in N$. If $P(h)=0$, then $ 0\neq Q(h)x\in N$. Therefore, $N$ contains a nonzero vector $v_3=R(h)x\in N$, where $\text{deg}_hR(h)=l$ and its leading coefficient‌ is $r_l\neq 0$. Furthermore, by \eqref{ae-a}, we have $0\neq (-2a)^ll!r_lx=(e-a)^lR(h)x\in N$. It follows that $x\in N$. Consequently, $\bar w_{\zeta}=\frac{1}{c}(E+s)x\in N$, which implies that $N=L(s,c,\xi_2)$.
\end{proof}
Next, we construct the corresponding simple Whittaker $\mathfrak g$-supermodules.
\begin{definition}
    Let $a,c\in \mathbb C^*$, and let $\mathbb C_{\chi}=\mathbb Cv_{\chi}$ be the 1-dimensional purely even \(\mathfrak l\)-supermodule defined by 
    \[
    ev_{\chi}=av_{\chi},\qquad\bar hv_{\chi}=cv_{\chi},\qquad \bar ev_{\chi}=\bar Ev_{\chi}=\bar Fv_{\chi}=0.
    \]
    Let $W(c)=\text{Ind}^{\mathrm U(\mathfrak g)}_{\mathrm U(\mathfrak l)}\mathbb C_\chi$, with an even cyclic vector \(w_0=1\otimes v_{\chi}\). As \(z_2\) is even and central, \((z_2-\xi_2)W(c)\) is a subsupermodule. Define $S(a,c,\xi_2)=W(c)/(z_2-\xi_2)W(c)$, and denote the image of \(w_0\) by $s_{\zeta}$. Then 
    $$
    es_{\zeta}=as_{\zeta},\qquad\bar es_{\zeta}=0,\qquad \bar hs_{\zeta}=cs_{\zeta},\qquad\bar Es_{\zeta }=\bar Fs_{\zeta}=0,\qquad z_2s_{\zeta }=\xi_2s_{\zeta }.
    $$
\end{definition}
    By Proposition \ref{center}, we have
    $$
    \begin{aligned}
        \xi_2s_{\zeta}=(2f\bar e+2\bar fe+h\bar h+\bar h+F\bar E+\bar FE)s_{\zeta}=2a\bar fs_{\zeta}+c(h+2)s_{\zeta}.
    \end{aligned}
    $$
    Thus, we have
    $$
    \bar fs_{\zeta}=\frac{\xi_2-c(h+2)}{2a}s_{\zeta}.
    $$
    Thus, \(\bar f\) can be eliminated from the PBW spanning set.
\begin{proposition}\label{LScong}
    Let \(a,c\in \mathbb C^*\) and choose $s\in \mathbb C^*$ such that $s^2=a$. Then $S(a,c,\xi_2)$ admits a PBW basis consisting of vectors
    $$
    \left\{f^mh^nF^\varepsilon E^\delta s_{\zeta}|m,n\in\mathbb{Z}_+,\varepsilon,\delta\in \{0,1\}\right\}.
    $$
Moreover, after forgetting the $\mathbb Z_2$-grading, $S(a,c,\xi_2)\cong L(s,c,\xi_2)$ as a \(\mathfrak g\)-module.
\end{proposition}
\begin{proof}
    By the PBW theorem, the supermodule \(W(c)\) has a vector space basis \[\{f^mh^nF^\varepsilon E^\delta\bar f^kw_0\mid m,n,k\in\mathbb{Z}_+, \varepsilon, \delta =0, 1\}\]
    The relation \(z_2w_0=2a\bar fw_0+c(h+2)w_0\), together with the centrality of \(z_2\), shows that $S(a,c,\xi_2)$ is spanned by the vectors $f^mh^nF^\varepsilon E^\delta s_{\zeta}$, where $m,n\geq 0$ and $\varepsilon ,\delta\in \{0,1\}$. It remains to prove that these vectors are linearly independent. Consider the simple $\mathfrak g$-module $L(s,c,\xi_2)$ with cyclic vector $\bar w_{\zeta}$ and set $x=\bar F\bar w_{\zeta }.$ Then we have
    \[
    ex=ax,\qquad\bar ex=\bar Ex=\bar Fx=0,\qquad\bar hx=cx,\qquad z_2x=\xi_2x.
    \]
    Thus, \(x\) satisfies the defining relations of \(s_\zeta\). The universal properties of induction and the central quotient yield a $\mathfrak{g}$-module homomorphism
    \[
    \theta:S(a,c,\xi_2)\to L(s,c,\xi_2)
    \]
    satisfying $\theta(s_{\zeta})=x$. Since $(E+s)x=c\bar w_{\zeta}$ and $c\neq 0$, the vector \(x\) is cyclic in \(L(s,c,\xi_2)\). Hence, \(\theta\) is surjective. Furthermore, by $Ex=c\bar w_{\zeta}-sx$, we have
    \begin{equation}\label{Theta}
        \theta\left(f^mh^nF^\varepsilon s_{\zeta}\right)=f^mh^nF^\varepsilon x,\qquad \theta\left(f^mh^nF^\varepsilon Es_{\zeta}\right)=f^mh^nF^\varepsilon(c\bar w_{\zeta}-sx).
    \end{equation}
    On the other hand, by Lemma \ref{Lbas}, we have
    \[
    \left\{f^mh^nF^\varepsilon \bar w_{\zeta},f^mh^nF^\varepsilon x|m,n\in\mathbb{Z}_+,\ \varepsilon\in \{0,1\}\right\}
    \] 
    forms a PBW basis of \(L(s,c,\xi_2)\). Note that 
    $$
    \begin{pmatrix}x\\Ex\end{pmatrix}=\begin{pmatrix}0&1\\c&-s\end{pmatrix}\begin{pmatrix}\bar w_{\zeta}\\x\end{pmatrix}.
    $$
    The determinant of the transformation matrix is $-c\neq0.$ Therefore, the vectors appearing on the right hand side of \eqref{Theta} are linearly independent, which implies that the vectors $f^mh^nF^\varepsilon E^\delta s_{\zeta}$, where $m,n\in\mathbb{Z}_+$ and $\varepsilon,\delta\in \{0,1\}$, are linearly independent. Consequently, these vectors form a basis of \(S(a,c,\xi_2)\), and \(\theta\) is an isomorphism of $\mathfrak g$-modules.
\end{proof}
\begin{theorem}\label{Ssim}
    For $a,c\in \mathbb C^*$ and $\xi_2\in\mathbb C$, the $\mathfrak g$-supermodule $S(a,c,\xi_2)$ is a simple Whittaker $\mathfrak g$-supermodule of type $\zeta$.
\end{theorem}
\begin{proof}
    Choose \(s\in\mathbb C^*\) with \(s^2=a\). By Proposition \ref{LScong} and Lemma \ref{sim2}, \( S(a,c,\xi_2)\) is simple as a $\mathfrak g$-module. Thus, $S(a,c,\xi_2)$ is a simple $\mathfrak g$-supermodule. Its even cyclic vector \(s_{\zeta}\) satisfies
    \[
    es_{\zeta}=as_{\zeta},\qquad
    \bar es_{\zeta}=0.
    \]
    Since \(s_{\zeta}\) is a cyclic Whittaker vector of type \(\zeta\), the assertion follows.
\end{proof}

We next establish that this construction exhausts the simple Whittaker $\mathfrak g$-supermodules with the prescribed central character.
\begin{lemma}\label{Vstru}
    Let \(V\) be a simple Whittaker $\mathfrak g$-supermodule of type $\zeta$ such that
    \[
    \zeta(e)=a\neq0,\qquad\zeta(\bar e)=0,\qquad z_1|_V=\xi_1\operatorname{id}_V, \qquad z_2|_V=\xi_2\operatorname{id}_V,\]
    where $\xi_1\neq0.$ Then there exists \(c\in\mathbb C^*\) with $c^2=\xi_1$ and a nonzero homogeneous cyclic Whittaker vector \(v\in V\) such that 
    \[
    ev=av,\qquad \bar hv=cv,\qquad \bar ev=\bar Ev=\bar Fv=0,\qquad z_2v=\xi_2v.
    \]
\end{lemma}
\begin{proof}
    Let \(v_0\) be a nonzero homogeneous Whittaker vector of type \(\zeta\) in \(V\). Since $\bar E^2=0$ and $[e,\bar E]=[\bar e,\bar E]=0$, Lemma \ref{cons} yields a nonzero homogeneous cyclic Whittaker vector \(v_1\) such that $\bar Ev_1=0$. Furthermore, by Proposition \ref{center}, we have $ z_1v_1=\bar h^2v_1=\xi_1v_1$. Choose $c_0\in \mathbb C^*$ such that $c_0^2=\xi_1.$ Then $(\bar h-c_0)(\bar h+c_0)v_1=0$. 
    Define
    \[
    v_+=\frac12\left(1+\frac{\bar h}{c_0}\right)v_1,\qquad v_-=\frac12\left(1-\frac{\bar h}{c_0}\right)v_1.
    \]
    Since $v_1=v_++v_-\neq 0$, at least one of \(v_+\) and \(v_-\) is nonzero. Moreover, 
    \[
    \bar hv_+=c_0v_+,\qquad \bar hv_-=-c_0v_-.
    \]
    Since $[e,\bar h]=-2\bar e$ and $[\bar e,\bar h]=[\bar E,\bar h]=0$, both $v_+$ and $v_-$, whenever nonzero, are Whittaker vectors of type $\zeta$ annihilated by \(\bar E\). Furthermore, since \(\bar h\) is even, both vectors are homogeneous. Consequently, there exists $c\in\{c_0,-c_0\}$ and a nonzero homogeneous Whittaker vector \(v_2\) such that
    \[
    \bar hv_2=cv_2,\qquad\bar Ev_2=0.
    \]
    Since $V$ is simple, \(v_2\) is also cyclic. Set $x=\bar Fv_2$. Then we have
    \[
    ex=ax,\qquad \bar ex=0,\qquad \bar Ex=\bar Fx=0,\qquad \bar hx=cx.
    \]
    If $x\neq 0$, it is a homogeneous Whittaker vector satisfying the required relations. Define
    \[
    v=
    \begin{cases}
        v_2,&\bar Fv_2=0,\\
        \bar Fv_2,&\bar Fv_2\neq0.
    \end{cases}
    \]
    Then $v$ is the required nonzero homogeneous cyclic Whittaker vector.
\end{proof}
We recall the parity reversal functor \(\Pi\) on \(\mathfrak g\)-supermodules. For a \(\mathfrak g\)-supermodule \(M=M_{\bar0}\oplus M_{\bar1}\), the parity-reversed supermodule \(\Pi M\) has the same underlying vector space as \(M\), with the grading
\[
(\Pi M)_{\bar i}=M_{\overline{i+1}},
\qquad i\in\{0,1\},
\]
and the \(\mathfrak g\)-action given by
\[
x\cdot\Pi(m)=(-1)^{|x|}\Pi(xm)
\]
for all homogeneous \(x\in\mathfrak g\) and \(m\in M\). More generally, for \(\varepsilon\in\{0,1\}\), we have 
\[
x\cdot\Pi^\varepsilon(m)
=
(-1)^{\varepsilon|x|}
\Pi^\varepsilon(xm).
\]
With this notation in place, we now proceed to the classification of simple Whittaker supermodules.
\begin{theorem}\label{theo1}
    Let $V$ be a simple Whittaker $\mathfrak g$-supermodule of type $\zeta$ satisfying 
    $$
        \zeta(e)=a\neq0,\qquad\zeta(\bar e)=0,\qquad z_1|_V=\xi_1\operatorname{id}_V, \qquad z_2|_V=\xi_2\operatorname{id}_V,
    $$
    where $\xi_1\neq0.$ Then there exists \(c\in\mathbb C^*\) and $\varepsilon\in \{0,1\}$ such that $c^2=\xi_1$ and $V\cong \Pi^{\varepsilon}S(a,c,\xi_2)$ as $\mathfrak g$-supermodules.
\end{theorem}
\begin{proof}
    By Lemma \ref{Vstru}, there exists a nonzero homogeneous cyclic Whittaker vector \(v\in V\) satisfying 
    \[
    ev=av,\qquad \bar hv=cv,\qquad \bar ev=\bar Ev=\bar Fv=0,\qquad z_2v=\xi_2v,
    \]
    where $c^2=\xi_1.$ Let $\varepsilon=|v|\in\{0,1\}.$ Define an even $\mathfrak g$-supermodule homomorphism   
    $$
    \pi:\Pi^\varepsilon S(a,c,\xi_2)\to V\qquad \Pi^{\varepsilon}(xs_{\zeta})\mapsto (-1)^{\varepsilon |x|}xv,\quad x\in\mathrm{U}(\mathfrak{g}).
    $$
    The defining relations of \(S(a,c,\xi_2)\) and the equality \(|v|=\varepsilon\) show that \(\pi\) is a well-defined even homomorphism. As \(v\) is cyclic, \(\pi\) is surjective. The simplicity of $\Pi^{\varepsilon}S(a,c,\xi_2)$, by Theorem \ref{Ssim}, implies $\ker\pi=0.$ Hence, $V\cong\Pi^\varepsilon S(a,c,\xi_2).$
\end{proof}
It remains to determine when two supermodules in this family are isomorphic.
\begin{theorem}\label{3.17}
    Let $a,c,d\in\mathbb C^*,\xi_2\in\mathbb C,$ and $\varepsilon,\delta\in\{0,1\}.$ Then $\Pi^\varepsilon S(a,c,\xi_2)\cong\Pi^\delta S(a,d,\xi_2)$ as \(\mathfrak g\)-supermodules if and only if $c=d$ and $\varepsilon=\delta.$
\end{theorem}
\begin{proof}
    Choose $s\in \mathbb C^*$ such that $s^2=a$. By an argument similar to that used in the proof of Lemma \ref{sim2}, we obtain that 
    $$
    \{v\in L(s,c,\xi_2)|(e-a)v=\bar ev=\bar Ev=\bar Fv=0\}=\mathbb C\bar F\bar w_{\zeta}.
    $$

    Furthermore, by Proposition \ref{LScong}, there exists a $\mathfrak g$-module isomorphism
    $$
    \theta :S(a,c,\xi_2)\to L(s,c,\xi_2)
    $$
    satisfying $\theta(s_{\zeta})=x$. Consequently, we have
    $$
    \{v\in S(a,c,\xi_2)|(e-a)v=\bar ev=\bar Ev=\bar Fv=0\}=\mathbb Cs_{\zeta}.
    $$

    Set $S=\Pi^{\varepsilon}S(a,c,\xi_2)$ and $S'=\Pi^{\delta}S(a,d,\xi_2)$. Let $v=\Pi^{\varepsilon}s_c$ and $v'=\Pi^{\delta}s_d$ be their respective homogeneous cyclic vectors, where $s_c$ and $s_d$ are the even cyclic vectors of $S(a,c,\xi_2)$ and $S(a,d,\xi_2)$. Then $|v|=|\varepsilon|,|v'|=|\delta|$.

    Suppose that there exists an even $\mathfrak g$-supermodule isomorphism $\tau:S\to S'.$ Since \(\tau\) is a homomorphism of \(\mathfrak g\)-supermodules, we have
    $$
    \begin{gathered}
        (e-a)\tau(v)=\bar e\tau(v)=\bar E\tau(v)=\bar F\tau(v)=0.
    \end{gathered}
    $$
    By an argument similar to that for $S(a,c,\xi_2)$ (Consider the actions of $\bar e$, $\bar E$, $\bar F$ and $e-a$ on $S'$ in sequence), we have
    $$
    \{v\in S'|(e-a)v=\bar ev=\bar Ev=\bar Fv=0\}=\mathbb Cv'.
    $$

    Therefore, there exists \(\eta\in\mathbb C^*\) such that 
    $\tau(v)=\eta v'$.
    Since $\bar hv=cv$, we have 
    $$
    \eta dv'=\eta \bar hv'=\bar h\tau(v)=\tau(\bar hv)=c\tau(v)=\eta cv'.
    $$
    Since \(\eta\ne0\) and \(v'\ne0\), it follows that $c=d$. Since \(\tau\) is even, $\delta=|v'|=|\tau(v)|=|v|=\varepsilon.$
    Thus, $\varepsilon=\delta$. 
    
    The sufficiency of the conditions is clear. Hence,
    $$
        \Pi^\varepsilon S(a,c,\xi_2)\cong\Pi^\delta S(a,d,\xi_2)\iff c=d,\ \varepsilon=\delta.
    $$
\end{proof}
\subsection{Whittaker supermodules with $b\neq 0$}\label{s33}
Assume $b\neq 0$, with $a$ arbitrary. Recall that the Universal Whittaker supermodule $M(\zeta)=\mathrm U(\mathfrak g)\otimes _{\mathrm U(\mathfrak n^+_{\bar 0})}\mathbb Cv_{\zeta}$, where $ev_{\zeta} = av_{\zeta}$ and $\bar ev_{\zeta} = bv_{\zeta}$. We determine \(\mathrm{Wh}_{\zeta}(M(\zeta))\) and construct an associative superalgebra acting on it. This yields a classification of the simple Whittaker $\mathfrak g$-supermodules of type \(\zeta\). 
\begin{proposition}\label{Mp}
    Define the following odd elements of $\mathrm U(\mathfrak g)$ by: 
    $$
    R=F-\frac{1}{2b}((h-\frac{a}{b}\bar h)\bar E+\bar hE),\qquad \overline{R}=\bar F-\frac{1}{2b}\bar h\bar E.
    $$ 
    Then the space of all Whittaker vectors in \(M(\zeta)\) is given by
    $$
    \mathrm{Wh}_{\zeta}(M(\zeta))=\bigoplus_{\substack{i,j\geq 0\\\alpha ,\beta ,\gamma,\delta\in \{0,1\}}}\mathbb Cz_1^iz_2^jE^{\alpha}\bar E^{\beta}R^{\gamma}\overline{R}^{\delta}v_{\zeta}.
    $$
\end{proposition}
\begin{proof}
    The commutation relations $[e,E]=[e,\bar E]=[\bar e,E]=[\bar e,\bar E]=0$ show that \(E\) and \(\bar E\) preserve Whittaker vectors. In particular, \(Ev_{\zeta}\) and \(\bar Ev_{\zeta}\) are Whittaker vectors. Furthermore, $[e,R]v_{\zeta}=[\bar e,R]v_{\zeta}=[e,\overline{R}]v_{\zeta}=[\bar e,\overline{R}]v_{\zeta}=0$ shows that \(R\) and \(\overline{R}\) also preserve Whittaker vectors. Since \(z_1\) and \(z_2\) are even and central, every vector of the form
    $$
    z_1^{i}z_2^{j}E^{\alpha }\bar E^{\beta }R^{\gamma}\overline{R}^{\delta}v_{\zeta}
    $$
    is a Whittaker vector. 
    Set
    \[
        \overline {\mathrm{Wh}_{\zeta}(M(\zeta))}=\operatorname{Span}_{\mathbb C}\left\{z_1^iz_2^jE^\alpha\bar E^\beta R^\gamma\overline{R}^\delta v_{\zeta}\right\},
    \]
    where \(i,j\geq0\) and \(\alpha,\beta,\gamma,\delta\in\{0,1\}\). Then $\overline{\mathrm{Wh}_{\zeta}(M(\zeta))} \subseteq \mathrm {Wh}_{\zeta}(M(\zeta))$. 
    
    To prove the reverse inclusion, observe that $\bar F\bar Ev_{\zeta}=\overline{R}\bar Ev_{\zeta}$ and $(R\bar E-E\overline{R})v_{\zeta}=(F\bar E-E\bar F)v_{\zeta}$. By the definitions of $R,\overline{R}$ and Proposition \ref{center}, we have 
    \begin{align}
    Rv_{\zeta} =& Fv_{\zeta}-\frac{1}{2b}\left( (h-\frac{a}{b}\bar{h})\bar{E}+\bar{h}E \right)v_{\zeta}, &  z_1v_{\zeta} =& (\bar{h}^2+2\overline{R}\bar{E}+4b\bar{f})v_{\zeta}, \label{R} \\[1ex]
    \overline{R}v_{\zeta} =& \bar{F}v_{\zeta}-\frac{1}{2b}\bar{h}\bar{E}v_{\zeta},
    &z_2v_{\zeta} =& (h\bar{h}+2bf+2\bar{h}+2a\bar{f}+R\bar{E}-E\overline{R})v_{\zeta}. \label{fc}
    \end{align}
    Then we have the decomposition which will be proved in Remark \ref{rem:PBWcom}:
    \begin{equation}\label{Mphi}
    M(\zeta)=\bigoplus\limits_{\substack{m,n,i,j\geq 0\\\alpha ,\beta ,\gamma,\delta\in \{0,1\}}}\mathbb Ch^m\bar h^nz_1^iz_2^jE^{\alpha}\bar E^{\beta}R^{\gamma}\overline{R}^{\delta}v_{\zeta}.    
    \end{equation}
    In particular, every vector $v\in M(\zeta)$ can be uniquely expressed as a finite sum
    $$
    v=\sum_{k}P_k(h,\bar h)u_k,\quad P_k(h,\bar h)\in \mathbb C[h,\bar h],\quad u_k\in \overline{\mathrm{Wh}_{\zeta}(M(\zeta))}.
    $$
    where the vectors $u_k$ are chosen from the fixed monomial basis in \eqref{Mphi}. Now suppose that $v\in \mathrm{Wh}_{\zeta}(M(\zeta))$. Then
    $$
    0=(\bar e-b)v=b\sum_{k}(P_k(h-2,\bar h)-P_k(h,\bar h))u_k.
    $$
    Since $b\neq 0$, uniqueness of the expansion \eqref{Mphi} gives $P_{k}(h-2,\bar h)=P_k(h,\bar h)$ for every $k$. Thus, \(P_k(h,\bar h)\) is independent of \(h\) and can be written as $P_k(\bar h)$.
    
    Similarly, we have 
    $$
    0=( e-a)v=-2b\sum_{k}\frac{\partial P_k}{\partial \bar h}u_k.
    $$
    Since $b\neq 0$ and the $\{u_{k}\}$ are linearly independent, $\frac{\partial P_k}{\partial \bar h}=0$ for every $k$. Hence, every $P_k(\bar h)$ is constant, proving the reverse inclusion. Consequently, $\mathrm{Wh}_{\zeta}(M(\zeta))=\overline{\mathrm{Wh}_{\zeta}(M(\zeta))}$.
\end{proof}
\begin{remark}\label{rem:PBWcom}
    We establish \eqref{Mphi} using a weighted filtration on \(\mathrm U(\mathfrak g)\) by
    \[
    \begin{gathered}
    \deg e=\deg\bar e=0,\qquad
    \deg h=\deg\bar h=1,\\
    \deg E=\deg\bar E=1,\qquad
    \deg F=\deg\bar F=3,\\
    \deg f=\deg\bar f=5.
    \end{gathered}
    \]

    Since every nonzero supercommutator of two generators has degree strictly smaller than the sum of their degrees, the PBW theorem yields that
    \[
    \operatorname{gr}M(\zeta)\cong S(f,\bar f,h,\bar h)\otimes\Lambda(E,\bar E,F,\bar F)
    \]
    as graded vector superspaces. Here \(S(f,\bar f,h,\bar h)\) denotes the symmetric algebra on the vector space \(\operatorname{span}_{\mathbb C}\{f,\bar f,h,\bar h\}\), and \(\Lambda(E,\bar E,F,\bar F)\) denotes the exterior algebra on \(\operatorname{span}_{\mathbb C}\{E,\bar E,F,\bar F\}\).

    By \eqref{R} and \eqref{fc}, the leading terms of $z_1v_{\zeta},z_2v_{\zeta},Rv_{\zeta},\bar Rv_{\zeta}$ are $4b\bar f,2bf+2a\bar f,F,\bar F$, respectively.
    Since \(b\neq0\), the linear transformation
    \[
    (f,\bar f)\mapsto(4b\bar f,\;2bf+2a\bar f)
    \]
    is invertible, with determinant \(-8b^2\neq0\). Thus, the leading terms of 
    \[
        h^m\bar h^n z_1^iz_2^jE^\alpha\bar E^\beta R^\gamma\bar R^\delta v_\zeta
    \]
    are
    \[
        h^m\bar h^n(4b\bar f)^i(2bf+2a\bar f)^jE^\alpha\bar E^\beta F^\gamma\bar F^\delta v_{\zeta},
    \]
    which form a basis of \(\operatorname{gr}M(\zeta)\), where $m,n,i,j\geq 0$ and $\alpha,\beta,\gamma,\delta\in \{0,1\}$. By the standard filtered-basis argument, these vectors form a basis of \(M(\zeta)\). Therefore,
\[
M(\zeta)=
\bigoplus\limits_{\substack{m,n,i,j\geq0\\
\alpha,\beta,\gamma,\delta\in\{0,1\}}}
\mathbb C h^m\bar h^n z_1^iz_2^j
E^\alpha\bar E^\beta
R^\gamma\bar R^\delta v_\zeta.
\]
\end{remark}
We next use an associated Clifford-type superalgebra to describe the simple Whittaker $\mathfrak g$-supermodules.
\begin{definition}
    Assume that \(b\neq0\). Define the following odd elements of $\mathrm U(\mathfrak g)$:
    $$
    p=\frac{1}{2b}\bar E,\quad r=E-\frac{a}{2b}\bar E,\quad s=R-\frac{1}{2b}\bar E,\quad q=\overline{R}.
    $$
    For $v\in \mathrm{Wh}_{\zeta}(V)$, we have 
    \[
    \begin{aligned}
    s^2v&=\left(
-f-\frac{F\bar E-E\bar F}{2b}
-\frac{(h+2)\bar h}{2b}
+\frac{a}{4b^2}\bar h^2
+\frac{a}{2b^2}\bar F\bar E
\right)v\\
&=\left(\frac{a}{4b^2}z_1-\frac1{2b}z_2\right)v,\\
[s,q]v
&=-\left(
2\bar f+\frac1b\bar F\bar E
+\frac1{2b}\bar h^2
\right)v\\
&=-\frac1{2b}z_1v.
\end{aligned}
\]
    Hence, their induced actions on Whittaker vectors satisfy the relations:
    \begin{equation}\label{rela1}
\begin{gathered}
    p^2=r^2=q^2=0, \quad s^2=\frac{a}{4b^2}z_1-\frac{1}{2b}z_2, \\
    [p, r]=1, \qquad [q, s]=-\frac{1}{2b}z_1, \\
    [p, s]=[p, q]=[r, s]=[r, q]=0.
\end{gathered}
\end{equation}
Let \(\mathcal A_\zeta\) denote the associative superalgebra generated by the even central elements \(z_1,z_2\) and the odd elements \(p,r,s,q\), subject to the defining relations \eqref{rela1}. 
\end{definition}
\begin{remark}
    The superalgebra $\mathcal A_{\zeta}$ is a Clifford superalgebra over $\mathbb C[z_1,z_2]$. More precisely, for a prescribed central character $(\xi_1,\xi_2)\in \mathbb C^2$, consider the central reduction 
    $$
    \mathcal C_{\xi_1,\xi_2}=\mathcal A_{\zeta}/\langle z_1-\xi_1,z_2-\xi_2\rangle.
    $$
    Write $c=\frac{a}{4b^2}\xi_1-\frac{1}{2b}\xi_2$ and $d=-\frac{1}{2b}\xi_1$. Then \(\mathcal C_{\xi_1,\xi_2}\) is the Clifford superalgebra associated with the four-dimensional purely odd vector space $U=\operatorname{Span}_{\mathbb C}\{p,r,s,q\}$ equipped with the symmetric bilinear form \(B_{\xi_1,\xi_2}\) determined by
    \[
    B_{\xi_1,\xi_2}(x,y)=xy+yx,\qquad x,y\in U.
    \]
    With respect to the ordered basis \((p,r,s,q)\), its matrix is
    \[
    B_{\xi_1,\xi_2}=
    \begin{pmatrix}
    0&1&0&0\\
    1&0&0&0\\
    0&0&2c&d\\
    0&0&d&0
    \end{pmatrix}.
    \]
    Since $\det B_{\xi_1,\xi_2}=d^2=\frac{\xi_1^2}{4b^2}$, the form \(B_{\xi_1,\xi_2}\) is nondegenerate if and only if \(\xi_1\neq0\). If \(\xi_1=0\), its radical is given by
    \[
    \operatorname{Rad}(B_{\xi_1,\xi_2})=
    \begin{cases}
        \mathbb Cq,&\xi_2\neq0,\\
        \mathbb Cs\oplus\mathbb Cq,&\xi_2=0.
    \end{cases}
    \]
    In the degenerate cases, the ideal generated by \(\operatorname{Rad}(B_{\xi_1,\xi_2})\) is nilpotent and therefore annihilates every simple \(\mathcal C_{\xi_1,\xi_2}\)-supermodule. Consequently, the classification of simple \(\mathcal C_{\xi_1,\xi_2}\)-supermodules reduces to that of the nondegenerate Clifford superalgebra associated with \(U/\operatorname{Rad}(B_{\xi_1,\xi_2})\).
\end{remark}
\begin{proposition}\label{Aphi}
    The superalgebra $\mathcal A_{\zeta}$ admits the following decomposition:
    $$
    \mathcal A_{\zeta}=\bigoplus\limits_{\alpha ,\beta ,\gamma,\delta\in \{0,1\}}\mathbb C[z_1,z_2]p^{\alpha}r^{\beta}s^{\gamma}q^{\delta}.
    $$ 
    Moreover, the map
    \[
        \mathcal A_\zeta\to \mathrm{Wh}_{\zeta}(M(\zeta)),
        \qquad u\mapsto uv_{\zeta},
    \]
    is an isomorphism of vector superspaces.
\end{proposition}
\begin{proof}
    By the definitions of $p,r,s,q$, we have
    $$
    \bar E=2bp,\qquad E=r+ap,\qquad R=s+p,\qquad \overline{R}=q.
    $$
    These formulas define an invertible triangular change of the odd generators. By Proposition \ref{Mp}, the vectors 
    \[
    p^\alpha r^\beta s^\gamma q^\delta v_{\zeta},\qquad \alpha,\beta,\gamma,\delta\in\{0,1\},
    \] 
    form a basis of $\mathrm{Wh}_{\zeta}(M(\zeta))$ over $\mathbb C[z_1,z_2]$. The relations \eqref{rela1} also show that these monomials span \(\mathcal A_\zeta\) over \(\mathbb C[z_1,z_2]\). Moreover, the ordered monomials \(p^\alpha r^\beta s^\gamma q^\delta\), where \(\alpha,\beta,\gamma,\delta\in\{0,1\}\), are linearly independent over \(\mathbb C[z_1,z_2]\), since the vectors \(p^\alpha r^\beta s^\gamma q^\delta v_\zeta\) form a basis of \(\mathrm{Wh}_{\zeta}(M(\zeta))\) over \(\mathbb C[z_1,z_2]\). Consequently, 
    $$
    \mathcal A_{\zeta}=\bigoplus\limits_{\alpha ,\beta ,\gamma,\delta\in \{0,1\}}\mathbb C[z_1,z_2]p^{\alpha}r^{\beta}s^{\gamma}q^{\delta}.
    $$
    
    The map \(u\mapsto uv_{\zeta}\) is therefore a parity-preserving linear bijection, and by Proposition \ref{Mp}, it is an isomorphism of vector superspaces.
\end{proof}
\begin{remark}\label{free}
    The induced module \(M(\zeta)\) admits a natural right \(\mathcal A_\zeta\)-module structure defined by
    \[
    (xv_{\zeta})\cdot u=xuv_{\zeta},\qquad x\in \mathrm U(\mathfrak g),\quad u\in\mathcal A_\zeta.
    \]
    Since \(uv_\zeta\) is a Whittaker vector of type \(\zeta\) for every \(u\in \mathcal A_\zeta\), we have $(e-a)uv_\zeta=(\bar e-b)uv_\zeta=0$. Thus, the above action respects the defining relations of the induced module \(M(\zeta)\) and is therefore well-defined.

    By Propositions \ref{Mp} and \ref{Aphi}, there is an isomorphism of right $\mathcal A_{\zeta}$-supermodules 
    \begin{equation*}
    M(\zeta)\cong \mathbb C[h,\bar h]v_{\zeta}\otimes \mathcal A_{\zeta}.
    \end{equation*}
    Thus, \(M(\zeta)\) is a free right \(\mathcal A_\zeta\)-module with basis $\{h^m\bar h^nv_{\zeta}\mid m,n\geq0\}.$
    This allows us to construct $\mathfrak g$-supermodules from \(\mathcal A_\zeta\)-supermodules.
\end{remark}
    For a left $\mathcal A_{\zeta}$-supermodule $S$, set $M(\zeta,S)=M(\zeta)\otimes _{\mathcal A_{\zeta}}S$, equipped with the left $\mathfrak g$-action $x\cdot (w_0\otimes s_0)=(xw_0)\otimes s_0$, for $x\in \mathrm U(\mathfrak g),w_0\in M(\zeta)$, and $s_0\in S.$ We first establish two properties of $M(\zeta,S)$ that will be used in the subsequent classification.
\begin{lemma}\label{forfil}
     There exists a natural isomorphism of vector superspaces $M(\zeta,S)\cong \mathbb C[h,\bar h]\otimes S.$ Moreover, the space of all Whittaker vectors in \(M(\zeta,S)\) is given by $\mathrm{Wh}_{\zeta}(M(\zeta,S))=v_{\zeta}\otimes S$.
\end{lemma}
\begin{proof}
    By Remark \ref{free}, $M(\zeta)$ is a free right $\mathcal A_{\zeta}$-module. Therefore, we obtain a natural isomorphism of vector superspaces
    $$
    M(\zeta,S)=M(\zeta)\otimes _{\mathcal A_{\zeta}}S\cong \mathbb C[h,\bar h]v_{\zeta}\otimes \mathcal A_{\zeta}\otimes_{\mathcal A_{\zeta}} S=\mathbb C[h,\bar h]v_{\zeta}\otimes S\cong \mathbb C[h,\bar h]\otimes S.
    $$
    Every vector in $v_{\zeta}\otimes S$ is a Whittaker vector. Conversely, by the argument similar to that used in the proof of Proposition \ref{Mp}, each $v\in \mathrm{Wh}_{\zeta}(M(\zeta,S))$ belongs to $v_{\zeta}\otimes S$. Hence, $\mathrm{Wh}_{\zeta}(M(\zeta,S))=v_{\zeta}\otimes S$.
\end{proof}
\begin{proposition}\label{subm}
    Let \(S\) be a nonzero \(\mathcal A_\zeta\)-supermodule. Then the following statements hold:
    \begin{enumerate}
        \item[(1).] Every nonzero \(\mathfrak g\)-subsupermodule of \(M(\zeta,S)\) contains a nonzero Whittaker vector.
        \item[(2).] If \(S\) is a simple \(\mathcal A_\zeta\)-supermodule, then \(M(\zeta,S)\) is a simple Whittaker $\mathfrak g$-supermodule.
    \end{enumerate}
\end{proposition}
\begin{proof}

Let \(N\) be a nonzero \(\mathfrak g\)-subsupermodule of \(M(\zeta,S)\), and denote by $\mathrm{Wh}_{\zeta}(N)$ the space of all Whittaker vectors in $N$.

(1) By Lemma \ref{forfil}, we have $M(\zeta,S)\cong\mathbb C[h,\bar h]\otimes S$. For each \(k\geq0\), define
    \[
        F_kM(\zeta,S)=\operatorname{Span}_{\mathbb C}\left\{h^m\bar h^nv_{\zeta}\otimes s\mid m+n\leq k,\ s\in S\right\}.
    \]
    Then \(\{F_kM(\zeta,S)\}_{k\geq 0}\) defines a filtration on \(M(\zeta,S)\). For any nonzero element \(v\in M(\zeta,S)\), we define its filtration degree by
\[
\deg(v):=\min\{k\geq 0\mid v\in F_kM(\zeta,S)\}.
\]
    Moreover, we have 
    $$
    (e-a)F_kM(\zeta,S)\subseteq F_{k-1}M(\zeta,S),\quad (\bar e-b)F_kM(\zeta,S)\subseteq F_{k-1}M(\zeta,S),
    $$
    where $F_{-1}M(\zeta,S)=0$. Choose a nonzero homogeneous vector \(v\in N\) of minimal filtration degree. Both $(e-a)v$ and $(\bar e-b)v$ lie in $N$ and have strictly smaller degree unless they vanish. Minimality therefore gives
    $$
    (e-a)v=(\bar e-b)v=0.
    $$
    Hence, $0\neq v\in \mathrm{Wh}_{\zeta}(N)$. 
    
(2) Suppose that \(S\) is a simple \(\mathcal A_\zeta\)-supermodule. By part (1) and Lemma \ref{forfil}, we have $0\neq \mathrm{Wh}_{\zeta}(N)\subset v_{\zeta}\otimes S$. Since \(\mathrm{Wh}_{\zeta}(N)\) is stable under the action of \(\mathcal A_\zeta\), it is a nonzero \(\mathcal A_\zeta\)-subsupermodule of \(v_{\zeta}\otimes S\). Since \(v_{\zeta}\otimes S\cong S\) as \(\mathcal A_\zeta\)-supermodules and \(S\) is simple, \(\mathrm{Wh}_{\zeta}(N)=v_{\zeta}\otimes S\). As these vectors generate \(M(\zeta,S)\), we obtain \(N=M(\zeta,S)\), proving simplicity.
\end{proof}
\begin{theorem}\label{Vcong}
    Let $V$ be a simple Whittaker $\mathfrak g$-supermodule of type $\zeta$. Then \(\mathrm{Wh}_{\zeta}(V)\) is a simple \(\mathcal A_\zeta\)-supermodule, and $V\cong M(\zeta)\otimes_{\mathcal A_\zeta}\mathrm{Wh}_{\zeta}(V)$ as \(\mathfrak g\)-supermodules.
\end{theorem}
\begin{proof}
    Define a surjective homomorphism
    $$
    \Phi:M(\zeta)\otimes_{\mathcal A_{\zeta}}\mathrm{Wh}_{\zeta}(V)\to V,\qquad xv_{\zeta}\otimes v\mapsto xv,
    $$
    for \(x\in \mathrm U(\mathfrak g)\) and \(v\in \mathrm{Wh}_{\zeta}(V)\). This map is a well-defined even $\mathfrak g$-supermodule homomorphism.
    
    If $\ker\Phi$ were nonzero, Proposition \ref{subm} provides a nonzero Whittaker vector in it. By Lemma \ref{forfil}, however, every Whittaker vector belongs to \(v_{\zeta}\otimes \mathrm{Wh}_{\zeta}(V)\). Since \(\Phi(v_{\zeta}\otimes v)=v\), the map \(\Phi\) is injective on this subspace.
    This contradiction yields $\ker\Phi=0$ and $V\cong M(\zeta)\otimes _{\mathcal A_{\zeta}}\mathrm{Wh}_{\zeta}(V)$ as $\mathfrak g$-supermodules. 
    
    Let \(T\) be a nonzero \(\mathcal A_\zeta\)-subsupermodule of \(\mathrm{Wh}_{\zeta}(V)\). Since \(M(\zeta)\) is a free right \(\mathcal A_\zeta\)-module, the tensor product
    \(
        M(\zeta)\otimes_{\mathcal A_\zeta}T
    \)
    is a nonzero \(\mathfrak g\)-subsupermodule of $
        M(\zeta)\otimes_{\mathcal A_\zeta}\mathrm{Wh}_{\zeta}(V)\cong V$. By the simplicity of \(V\), we obtain
    \[
    M(\zeta)\otimes_{\mathcal A_\zeta}T=M(\zeta)\otimes_{\mathcal A_\zeta}\mathrm{Wh}_{\zeta}(V).
    \]
    Taking Whittaker vectors on both sides and applying Lemma \ref{forfil}, we have $v_{\zeta}\otimes T=v_{\zeta}\otimes \mathrm{Wh}_{\zeta}(V)$. Thus, $T=\mathrm{Wh}_{\zeta}(V)$, proving its simplicity as an \(\mathcal A_\zeta\)-supermodule. 
\end{proof}
We next determine the simple $\mathcal A_{\zeta}$-supermodule structure of $\mathrm{Wh}_{\zeta}(V)$ by means of Clifford superalgebras. By the standard classification of complex superalgebras (see \cite[Section 3.1 and Exercise 3.11]{CW2012}),
\[
\mathrm{Cl}_2\cong M(1|1),\qquad
\mathrm{Cl}_3\cong Q(2),\qquad
\mathrm{Cl}_4\cong M(2|2).
\]

Here \(M(m|n)\) and \(Q(n)\) denote the matrix and queer simple associative superalgebras, respectively. A simple supermodule over an algebra of type \(M\) is not isomorphic to its parity reversal, whereas a simple supermodule over an algebra of type \(Q\) is isomorphic to its parity reversal.

We now obtain the following classification.
\begin{theorem}\label{Vphi}
    Let \(V\) be a simple Whittaker \(\mathfrak g\)-supermodule of type \(\zeta\) on which \(z_1,z_2\) act by scalars \(\xi_1,\xi_2\in\mathbb C\), respectively. Then there exists a nonzero homogeneous vector $u\in \mathrm{Wh}_{\zeta}(V)$ such that 
    \[
    \mathrm{Wh}_{\zeta}(V)=
    \begin{cases}
    \operatorname{Span}_{\mathbb C}\{u,ru\},
    &\text{ for }\xi_1=\xi_2=0,\\
    \operatorname{Span}_{\mathbb C}\{u,ru,su,rsu\},&\text{ otherwise}.
    \end{cases}
    \]
    The above spanning sets form bases of \(\mathrm{Wh}_{\zeta}(V)\). The even-isomorphism classes of these simple \(\mathcal A_\zeta\)-supermodules are as follows:
    \begin{itemize}
        \item[(1).] If \(\xi_1\neq0\) or $\xi_1=\xi_2=0$, there are exactly two isomorphism classes, interchanged by parity reversal.
        \item[(2).] If \(\xi_1=0\) and \(\xi_2\neq0\), there is exactly one isomorphism class, and \(\mathrm{Wh}_{\zeta}(V)\cong\Pi\mathrm{Wh}_{\zeta}(V)\).
    \end{itemize}
\end{theorem}
\begin{proof}
    By Theorem \ref{Vcong}, \(\mathrm{Wh}_\zeta(V)\) is a simple \(\mathcal A_\zeta\)-supermodule. Since \(z_1\) and \(z_2\) act by \(\xi_1\) and \(\xi_2\), respectively, the action factors through \(\mathcal C_{\xi_1,\xi_2}\). Thus, \(\mathrm{Wh}_\zeta(V)\) is a simple \(\mathcal C_{\xi_1,\xi_2}\)-supermodule. We distinguish three cases according to the rank of \(B_{\xi_1,\xi_2}\).

\medskip

\noindent\textbf{Case 1: $\xi_1\neq 0$.}
In this case, \(d\neq0\), so \(B_{\xi_1,\xi_2}\) is nondegenerate of rank \(4\). Hence
\[
\mathcal C_{\xi_1,\xi_2}
\cong \mathrm{Cl}_4
\cong M(2|2).
\]

By the classification of simple supermodules over matrix superalgebras, there are exactly two simple \(\mathcal C_{\xi_1,\xi_2}\)-supermodules up to even isomorphism, interchanged by parity reversal, and both have superdimension \((2|2)\).

More explicitly, set $\widetilde t=s-\frac cdq,\widetilde z=\frac qd$. Then
\[
p^2=r^2=\widetilde t^2=\widetilde z^2=0,
\qquad
[p,r]=[\widetilde t,\widetilde z]=1,
\]
and all other brackets among these generators vanish. By the standard simple Clifford supermodule construction, there exists a nonzero homogeneous vector \(u\in \mathrm{Wh}_{\zeta}(V)\) such that $pu=qu=0$, and the corresponding Clifford basis is $\{u,ru,\widetilde tu,r\widetilde tu\}$. Since \(\widetilde tu=su\), we obtain
\[
\mathrm{Wh}_\zeta(V)
=
\operatorname{Span}_{\mathbb C}
\{u,ru,su,rsu\}.
\]

\medskip

\noindent\textbf{Case 2:
$\xi_1=0$ and $\xi_2\neq0$.}
Here \(d=0\) and \(c\neq0\), and $\operatorname{Rad}(B_{\xi_1,\xi_2})=\mathbb Cq$. 
The ideal generated by \(q\) is nilpotent and hence annihilates every simple supermodule. Consequently, 
\[
\mathcal C_{\xi_1,\xi_2}/(q)
\cong\mathrm{Cl}_3
\cong Q(2).
\]

The queer superalgebra \(Q(2)\) has a unique simple supermodule up to even isomorphism, of superdimension \((2|2)\), which is isomorphic to its parity reversal.

Since \(p^2=0\), we may choose a nonzero homogeneous vector \(u\in \mathrm{Wh}_{\zeta}(V)\) with \(pu=0\). The relations \(s^2=c\neq0\) and \([p,s]=0\) imply that \(u,su\) are linearly independent and are both annihilated by \(p\). Using \([p,r]=1\) and the superdimension of the simple supermodule, we obtain 
\[
\mathrm{Wh}_\zeta(V)
=
\operatorname{Span}_{\mathbb C}
\{u,ru,su,rsu\},
\qquad
\mathrm{Wh}_\zeta(V)\cong\Pi\mathrm{Wh}_\zeta(V).
\]

\medskip

\noindent\textbf{Case 3:
$\xi_1=0$ and $\xi_2=0$.}
Here \(c=d=0\), and $\operatorname{Rad}(B_{\xi_1,\xi_2})=\mathbb Cs\oplus\mathbb Cq$.

The ideal generated by \(s\) and \(q\) is nilpotent and hence annihilates every simple supermodule. Consequently,
\[
\mathcal C_{\xi_1,\xi_2}/(s,q)
\cong\mathrm{Cl}_2
\cong M(1|1).
\]

The matrix superalgebra \(M(1|1)\) has exactly two simple supermodules up to even isomorphism, interchanged by parity reversal, and both have superdimension \((1|1)\).

Choosing a nonzero homogeneous vector \(u\in \mathrm{Wh}_{\zeta}(V)\) with \(pu=0\), we obtain 
\[
\mathrm{Wh}_\zeta(V)
=
\operatorname{Span}_{\mathbb C}\{u,ru\}.
\]

This completes the classification.
\end{proof}
\begin{theorem}\label{Vphi2}
Let $V$ be a simple Whittaker $\mathfrak g$-supermodule of type $\zeta$, where
\[
\zeta(e)=a,\qquad \zeta(\bar e)=b\neq 0.
\]
Suppose $z_1,z_2$ act on $V$ by $\xi_1$ and $\xi_2$, respectively. Then $V$ has the following $\mathfrak g$-supermodule description:
\[
    V\cong M(\zeta)\otimes _{\mathcal A_{\zeta}}\mathrm{Wh}_{\zeta}(V)\cong 
    \begin{cases}
    M(\zeta)\otimes _{\mathcal A_{\zeta}} \operatorname{Span}_{\mathbb C}\{u,ru\},
    &\text{ for }\xi_1=\xi_2=0,\\
    M(\zeta)\otimes _{\mathcal A_{\zeta}} \operatorname{Span}_{\mathbb C}\{u,ru,su,rsu\},&\text{ otherwise }.
    \end{cases}
\]
Here \(u\in \mathrm{Wh}_{\zeta}(V)\) is nonzero and homogeneous, and the given spanning sets are bases of \(\mathrm{Wh}_{\zeta}(V)\).

Moreover, the simple $\mathfrak g$-supermodules \(V\) are classified, up to even isomorphism, as follows.
    \begin{itemize}
        \item[(1).] If \(\xi_1\neq0\) or $\xi_1=\xi_2=0$, there are exactly two isomorphism classes, interchanged by parity reversal.
        \item[(2).] If \(\xi_1=0\) and \(\xi_2\neq0\), there is exactly one isomorphism class, and \(V\cong\Pi V\).
    \end{itemize}
\end{theorem}
\begin{proof}
    This is an immediate consequence of Lemma \ref{forfil}, Proposition \ref{subm}, Theorems \ref{Vcong} and \ref{Vphi}.
\end{proof}
\subsection{Whittaker supermodules with $a=b=0$}\label{s34}
We now consider the zero Whittaker character \(\zeta=0\). Throughout this subsection, let \(V\) be a simple Whittaker $\mathfrak g$-supermodule on which \(z_1\) and \(z_2\) act by \(\xi_1\) and \(\xi_2\), respectively. We begin with a refinement of the cyclic Whittaker vector.
\begin{proposition}\label{426}
    There exists a nonzero homogeneous cyclic Whittaker vector $\widetilde w\in V$ satisfying
    $$
    e\widetilde w=\bar e\widetilde w=E\widetilde w=\bar E\widetilde w=0,\qquad \bar h\widetilde w=\mu\widetilde w,\qquad \mu(h+1)\widetilde w=\xi_2 \widetilde w,
    $$
    where $\mu \in \mathbb C$ satisfies $\mu^2=\xi_1$. In particular, if $\xi_1=0$, we have
    $$
    e\widetilde w=\bar e\widetilde w=E\widetilde w=\bar E\widetilde w=\bar h\widetilde w=\bar F\widetilde w=0.
    $$
\end{proposition}
\begin{proof}
    Similar to the argument of Lemma \ref{Vstru}, by Lemma \ref{cons} and $\bar EE=(2\bar e-E\bar E)$, there exists a nonzero homogeneous cyclic Whittaker vector $\widetilde w_0$ with $E\widetilde w_0=\bar E\widetilde w_0=0$ and $\bar h^2\widetilde w_0=\mu ^2\widetilde w_0$ for some $\mu\in \mathbb C$ satisfying $\mu^2=\xi_1$. 
    
    If $\xi_1=0$, by Lemma \ref{cons} and $E\bar h=\bar hE-\bar E$, there exists a nonzero homogeneous cyclic Whittaker vector $\widetilde w_1$ with $e\widetilde w_1=\bar e\widetilde w_1=E\widetilde w_1=\bar E\widetilde w_1=\bar h\widetilde w_1=0$. And we have $\bar F\widetilde w_1$ is Whittaker vector with
    $$
    \bar F^2\widetilde w_1=E\bar F\widetilde w_1=\bar E\bar F\widetilde w_1=\bar h\bar F\widetilde w_1=0.
    $$
    Use Lemma \ref{cons} again. Then there exists a nonzero homogeneous cyclic Whittaker vector $\widetilde w$ with $e\widetilde w=\bar e\widetilde w=E\widetilde w=\bar E\widetilde w=\bar h\widetilde w=\bar F\widetilde w=0$.

    If $\xi_1\neq 0$, the argument of Lemma \ref{Vstru} yields a nonzero homogeneous cyclic Whittaker vector $\widetilde w$ with $E\widetilde w=\bar E\widetilde w=0$ and $\bar h\widetilde w=\mu \widetilde w$.
    Proposition \ref{center} then gives 
\begin{align*}
    \xi_2\widetilde w=z_2\widetilde w=(h+1)\bar h\widetilde w=\mu(h+1)\widetilde w.
\end{align*}
\end{proof}
\begin{remark}\label{428}
    If $\xi_1\neq 0$, then $\mu\neq 0$, and Proposition \ref{426} gives $h\widetilde w=\lambda\widetilde w, \lambda=\xi_2\mu^{-1}-1$.
    Thus, $\widetilde w$ is a highest-weight vector under the Cartan subalgebra $\mathfrak h$, and $V$ is a simple highest-weight $\mathfrak g$-supermodule. 
    In particular, when $\xi_1=0$, Proposition \ref{426} imposes a further restriction on the central character. Indeed,
    $$
        \xi_2\widetilde w=z_2\widetilde w=(h+1)\bar h\widetilde w=0.
    $$ 
    Since \(\widetilde w\neq0\), it follows that \(\xi_2=0\). Consequently, no simple Whittaker \(\mathfrak g\)-supermodule of type \(\zeta=0\) exists when \(\xi_1=0\) and \(\xi_2\neq0\).
\end{remark}
We now give an explicit classification of the simple Whittaker supermodules arising in Remark \ref{428}. Let $\mathfrak b^+=\mathfrak h\oplus\mathfrak n^+$. For \(\lambda,\mu\in\mathbb C\), let \(\mathbb C_{\lambda,\mu}\) be the one-dimensional \(\mathfrak b^+\)-supermodule generated by an even vector \(w_{\lambda,\mu}\) with
\[
\mathfrak n^+w_{\lambda,\mu}=0,\qquad
hw_{\lambda,\mu}=\lambda w_{\lambda,\mu},\qquad
\bar h w_{\lambda,\mu}=\mu w_{\lambda,\mu}.
\]
Define the Verma module
\[
\mathcal{M}(\lambda,\mu)
=\operatorname{Ind}_{\mathrm{U}(\mathfrak b^+)}^{\mathrm{U}(\mathfrak g)}
\mathbb C_{\lambda,\mu}.
\]
And we denote the simple quotient of $\mathcal{M}(\lambda,\mu)$ by $\mathcal{L}(\lambda,\mu)$. 
\begin{theorem}\label{high}
    Let \(V\) be a simple Whittaker \(\mathfrak g\)-supermodule of type \(\zeta=0\), on which \(z_1\) and \(z_2\) act by scalars \(\xi_1,\xi_2\in\mathbb C\), respectively. Suppose that \(\xi_1\neq0\). Then there exist \(\lambda,\mu\in\mathbb C\) and \(\varepsilon\in\{0,1\}\), with \(\mu^2=\xi_1,\lambda=\xi_2\mu^{-1}-1\), such that
    \[
    V\cong\Pi^\varepsilon \mathcal{L}\left(\lambda,\mu\right).
    \]

    Moreover, for \(\mu,\mu'\neq0\),
    \[
    \Pi^\varepsilon \mathcal{L}(\lambda,\mu)\cong\Pi^\delta \mathcal{L}(\lambda',\mu')\iff
    (\lambda,\mu)=(\lambda',\mu')
    \quad\text{and}\quad
    \varepsilon=\delta.
    \]

    In particular, for each fixed pair \((\xi_1,\xi_2)\) with \(\xi_1\neq0\), there are exactly four even isomorphism classes of simple Whittaker \(\mathfrak g\)-supermodules of type \(\zeta=0\).
\end{theorem}
\begin{proof}
    We first prove $\mathcal{M}(\lambda,\mu)=\mathcal{L}(\lambda,\mu)$ for $\mu\neq 0$.  By the PBW theorem and the defining relations, \(\bar h-\mu\) acts locally nilpotently on \(\mathcal{M}(\lambda,\mu)\). Let $N$ be a nonzero proper $\mathfrak g$-subsupermodule of \(\mathcal{M}(\lambda,\mu)\). Choose a nonzero homogeneous vector \(v_{\lambda,\mu}\in N\) of maximal \(h\)-weight \(\lambda-n\), where \(n>0\). Then $\mathfrak n^+v_{\lambda,\mu}=0$.
    Proposition \ref{center} then gives 
    $$
    \mu(\lambda+1)v_{\lambda,\mu}=z_2v_{\lambda,\mu}=(\lambda-n+1)\bar hv_{\lambda,\mu}.
    $$
    Consequently, 
    $$
    (\lambda-n+1)(\bar h-\mu)v_{\lambda,\mu}=n\mu v_{\lambda,\mu}.
    $$
    Since \(\bar h-\mu\) acts locally nilpotently on \(\mathcal{M}(\lambda,\mu)\), there exists a positive integer \(k\) such that $(\bar h-\mu)^k v_{\lambda,\mu}=0$. Hence, $0=(n\mu)^k v_{\lambda,\mu}$, which contradicts \(n\geq 1\), \(\mu\neq 0\) and \(v_{\lambda,\mu}\neq 0\). Hence, \(\mathcal{M}(\lambda,\mu)=\mathcal L(\lambda,\mu)\) for \(\mu\neq0\).

    By Proposition \ref{426}, \(V\) contains a nonzero homogeneous cyclic vector \(\widetilde w\) satisfying
    \[
    \mathfrak n^+\widetilde w=0,\qquad\bar h\widetilde w=\mu\widetilde w,\qquad h\widetilde w=\lambda\widetilde w,
    \]
    where \(\mu^2=\xi_1\neq0,\lambda=\xi_2\mu^{-1}-1\). Let \(\varepsilon\) be the parity of \(\widetilde w\). By the universal property of induced supermodules, there exists a surjective even \(\mathfrak g\)-supermodule homomorphism
    \[
    \Pi^\varepsilon
    \mathcal M(\lambda,\mu)\to  V.
    \]
    By the simplicity of $\mathcal{M}(\lambda,\mu)$, this homomorphism is an isomorphism. Hence, $V\cong \Pi^{\varepsilon}\mathcal M(\lambda,\mu)=\Pi^{\varepsilon}\mathcal L(\lambda,\mu)$.

    Finally, each \(\mathcal{L}(\lambda,\mu)\) with \(\mu\neq0\) is a simple Whittaker supermodule of type \(\zeta=0\), and Proposition \ref{center} gives its central character as
    $(\xi_1,\xi_2)=(\mu^2,\mu(\lambda+1))$.
    Since its highest-weight space is one-dimensional, an even isomorphism must preserve both the highest weight and the parity of the highest-weight vector. Thus, 
    \[
    \Pi^\varepsilon \mathcal{L}(\lambda,\mu)\cong\Pi^\delta \mathcal{L}(\lambda',\mu')\iff
    (\lambda,\mu)=(\lambda',\mu')
    \quad\text{and}\quad
    \varepsilon=\delta.
    \]
\end{proof}
Throughout the remainder of this subsection, we assume $\zeta=\xi_1=\xi_2=0.$ We describe \(V\) via modules over a 2-dimensional solvable Lie algebra.
\begin{definition}
    For any $\mathfrak g$-supermodule $M$, we define 
    $$
        \mathcal N(M)=\{v\in M|Ev=\bar Ev=\bar hv=\bar Fv=0\}.
    $$
\end{definition}
Proposition \ref{426} gives $\mathcal N(V)\neq 0$. The defining relations imply $e\mathcal N(V)=\bar e\mathcal N(V)=0$. Hence, every nonzero homogeneous vector of $\mathcal N(V)$ is a cyclic Whittaker vector by the simplicity of $V$.
The defining relations show that \(\mathcal N(V)\) is invariant under the action of \(h\) and \(\bar f\). Accordingly, consider the 2-dimensional nonabelian and solvable Lie subalgebra $\mathfrak s=\mathbb Ch\oplus\mathbb C\bar f.$ Then \(\mathcal N(V)\) is naturally an \(\mathfrak s\)-module. We consider the two cases of whether the action of \(\bar f\) on \(\mathcal N(V)\) is injective.
If \(\bar f\) has a nonzero kernel on \(\mathcal N(V)\), choose a nonzero homogeneous \(v\in \mathcal N(V)\) with \(\bar fv=0\). Since \(\mathfrak p\) is an ideal of \(\mathfrak g\), $\{v\in V|\mathfrak pv=0\}$ is nonzero subsupermodule. By the simplicity of $V$, it follows that \(\mathfrak pV=0\).
\begin{theorem}\label{0osp}
    Let $V$ be a simple Whittaker $\mathfrak g$-supermodule of type $\zeta=0$, on which $z_1$ and $z_2$ act as zero, and $\mathfrak pV=0$. Then there exist $\lambda\in \mathbb C$ and $\varepsilon\in \{0,1\}$ such that
    $$
    V\cong \Pi^{\varepsilon}\mathcal L(\lambda).
    $$
    Here $\mathcal L(\lambda)$ is the unique simple quotient supermodule of the Verma supermodule on $\mathfrak{osp}(1|2)$.
    Moreover,
    \[
    \Pi^\varepsilon \mathcal L(\lambda)\cong\Pi^\delta \mathcal L(\eta)\iff\lambda=\eta,\quad \varepsilon=\delta.
    \]
    Conversely, every \(\Pi^\varepsilon \mathcal L(\lambda)\) is a simple Whittaker supermodule of type \(\zeta=0\).
\end{theorem}
\begin{proof}
    Since $E^2=e$, by Lemma \ref{cons}, there exists a nonzero homogeneous Whittaker vector \(u\in V\) such that \(Eu=0\). Consider the Scasimir element $\Sigma=EF-FE+\frac{1}{2}$. Since \(\Sigma^2\) is central, it acts on \(V\) by a scalar \(\chi\in\mathbb C\). Moreover, \(Eu=0\) implies $\Sigma u=\left(h+\frac12\right)u$. Hence,
    \[
        ((h+\frac12)^2-\chi)u=0.
    \]
    Thus, \(\mathbb C[h]u\) contains a nonzero homogeneous vector \(v\) satisfying \(hv=\lambda v\) for some \(\lambda\in\mathbb C\). Since $Eh=(h-1)E$, we have \(Ev=0\). Hence, the vector \(v\) is a highest-weight vector. By simplicity, \(V\) is a simple highest-weight supermodule, so
    \[
        V\cong\Pi^\varepsilon \mathcal L(\lambda)
    \]
    for some \(\varepsilon\in\{0,1\}\).

    Conversely, every \(\Pi^\varepsilon \mathcal L(\lambda)\) is a simple Whittaker supermodule of type \(\zeta=0\). Finally, the highest weight and the parity of its highest-weight vector uniquely determine the even isomorphism class. This completes the proof.
\end{proof}
It remains to treat the case in which \(\bar f\) acts injectively on \(\mathcal N(V)\). Choose a nonzero homogeneous element \(u\in \mathcal N(V)\) and set \(S=\mathrm U(\mathfrak s)u\subseteq \mathcal N(V)\). Then \(\bar f\) acts injectively on \(S\).
Since \(\mathfrak s\) is purely even, all nonzero homogeneous elements of \(S\) have the same parity as \(u\). After forgetting this parity, we regard \(S\) as a purely even \(\mathfrak s\)-supermodule.
We extend the $\mathfrak s$-module structure on $S$ to an $\mathfrak a$-supermodule by requiring
$$
eS=ES=\bar eS=\bar ES=\bar hS=\bar FS=0.
$$
The defining relations of \(\mathfrak a\) show that this extension is well-defined. Define the induced \(\mathfrak g\)-supermodule
$$
   \mathrm{Ind}_{\mathfrak a}^{\mathfrak g}S:=\mathrm U(\mathfrak g)\otimes_{\mathrm U(\mathfrak a)}S.
$$ 
The PBW theorem gives the following decomposition of vector superspaces, with S identified with $1\otimes S$:
$$
\mathrm{Ind}_{\mathfrak a}^{\mathfrak g}S=\bigoplus\limits_{i\geq0}(f^iS\oplus f^iFS).
$$
We record three properties of this construction.
\begin{lemma}\label{nn}
    Assume that $\bar f$ acts injectively on $S$. Let $N$ be any nonzero $\mathfrak g$-subsupermodule of $\mathrm{Ind}_{\mathfrak a}^{\mathfrak g}S$. Then $N\cap S\neq0$.
\end{lemma}
\begin{proof}
    Let $v_N\in N$ be a nonzero vector. Its unique PBW expansion is $v_N=\sum\limits_{m=0}^{d_2}f^m(s_m+Ft_m)$ where $s_m,t_m\in S$. Suppose that $t_m\neq0$ for some $m$. Then 
    $$
        \bar Fv_N=-2\sum_{m=0}^{d_2}f^m\bar ft_m.
    $$
    The injectivity of $\bar f$ implies $0\neq\bar Fv_N\in N$. Replacing $v_N$ by $\bar Fv_N$ when necessary, we assume  
    $$
        v_N=\sum_{m=0}^{d_2}f^ms_m,\quad s_m\in S,s_{d_2}\neq0.
    $$
    If $d_2=0$, then $0\neq v_N=s_0\in N\cap S$. If $d_2>0$, the defining relations give $\bar h^{d_2}v_N=(-2)^{d_2}d_2!\bar f^{d_2}s_{d_2}$. Since $\bar f$ is injective and $s_{d_2}\neq 0$, there exists a nonzero element in $N\cap S$. Thus, $N\cap S\neq 0$.
\end{proof}
\begin{lemma}\label{eeqq}
    Under the assumptions of Lemma \ref{nn}, we have 
    $$
        \mathcal N(\mathrm{Ind}_{\mathfrak a}^{\mathfrak g}S)=S.
    $$
\end{lemma}
\begin{proof}
    The defining $\mathfrak a$-action gives $S\subseteq \mathcal N(\mathrm{Ind}_{\mathfrak a}^{\mathfrak g}S)$. We prove the reverse inclusion.

    Let $v\in \mathcal N(\mathrm{Ind}_{\mathfrak a}^{\mathfrak g}S)$. By the PBW theorem, we write $v=\sum\limits_{m=0}^{d_3}f^m(s_{m,0}+Fs_{m,1})$, where \(s_{m,0},s_{m,1}\in S\). Since $\bar Fv=0$, we obtain
    $$
    0=-2\sum_{m=0}^{d_3}f^m\bar fs_{m,1}.
    $$
    The injectivity of \(\bar f\) and linear independence of PBW vectors give $s_{m,1}=0$ for all $0\le m\le d_3$. Hence, $v=\sum\limits_{m=0}^{d_3}f^ms_{m,0}$.
    On the other hand, \(\bar hv=0\). Using \([\bar h,f]=-2\bar f\) and \(\bar hS=0\), we obtain
    $$
    0=-2\sum_{m=1}^{d_3}mf^{m-1}\bar fs_{m,0}.
    $$
    The same argument gives $s_{m,0}=0$ for $1\le m\le d_3$. Hence, $v=s_{0,0}\in S$, proving \(\mathcal N(\mathrm{Ind}_{\mathfrak a}^{\mathfrak g}S)=S\). 
\end{proof}

\begin{lemma}\label{simeq}
    If $\bar f$ acts on \(S\) injectively, then $\mathrm{Ind}_{\mathfrak a}^{\mathfrak g}S$ is a simple $\mathfrak g$-supermodule if and only if $S$ is a simple $\mathfrak s$-module.
\end{lemma}
\begin{proof}
    Suppose first that $S$ is a simple $\mathfrak s$-module and let $N$ be a nonzero $\mathfrak g$-subsupermodule of $\mathrm{Ind}_{\mathfrak a}^{\mathfrak g}S$.
    By Lemma \ref{nn}, we have $N\cap S\neq0.$
    Since both $N$ and $S$ are stable under the action of $\mathfrak s$, the intersection $N\cap S$ is a nonzero $\mathfrak s$-submodule of $S$. By the simplicity of $S$, we have $N\cap S=S.$
    Thus, $S\subseteq N$. Since $\mathrm{Ind}_{\mathfrak a}^{\mathfrak g}S=\mathrm U(\mathfrak g)S$, we conclude that $N=\mathrm{Ind}_{\mathfrak a}^{\mathfrak g}S$. Hence, $\mathrm{Ind}_{\mathfrak a}^{\mathfrak g}S$ is a simple $\mathfrak g$-supermodule.

    Conversely, suppose that $\mathrm{Ind}_{\mathfrak a}^{\mathfrak g}S$ is simple, and let $0\neq S'\subseteq S$ be an $\mathfrak s$-submodule. We regard $S'$ as an $\mathfrak a$-module by restricting the $\mathfrak a$-action on $S$, and consider the induced supermodule $\mathrm{Ind}_{\mathfrak a}^{\mathfrak g}S'$.
    Inclusion $S'\hookrightarrow S$ induces an injective homomorphism of $\mathfrak g$ -supermodule $\tau :\mathrm{Ind}_{\mathfrak a}^{\mathfrak g}S'\to \mathrm{Ind}_{\mathfrak a}^{\mathfrak g}S.$
    Its image is nonzero, therefore $\tau $ is surjective by the simplicity of \(\mathrm{Ind}_{\mathfrak a}^{\mathfrak g}S\). The PBW decompositions then imply \(S'=S\), proving the simplicity of $S$.
\end{proof}
We obtain the following structure theorem for \(\zeta=\xi_1=0\).
\begin{theorem}\label{cc}
    Let $V$ be a simple Whittaker $\mathfrak g$-supermodule of type $\zeta=0$, and suppose that $z_1,z_2$ act on V by zero. Then the following statements hold:
    \begin{itemize}
        \item[(1).] If $\mathfrak pV=0$, then there exist $\lambda\in \mathbb C$ and $\varepsilon\in \{0,1\}$ such that $V\cong \Pi^{\varepsilon}\mathcal L(\lambda)$.
        \item[(2).] If $\mathfrak pV\neq 0$, then there exists an element $\varepsilon\in \{0,1\}$ such that $V\cong \Pi^{\varepsilon}(\mathrm{Ind}_{\mathfrak a}^{\mathfrak g}S)$ as $\mathfrak g$-supermodules.
    \end{itemize}
\end{theorem}
\begin{proof}
    If $\mathfrak pV=0$, it is proved by Theorem \ref{0osp}.

    Now, we consider $\mathfrak pV\neq 0$. By the argument preceding Lemma \ref{nn}, \(\bar f\) acts injectively on \(\mathcal N(V)\). Set $\varepsilon=|u|$ and by the prescribed $\mathfrak a$-action on $S$ and the universal property of induced supermodule, there exists an even $\mathfrak g$-supermodule homomorphism
    $$
    \pi:\Pi^{\varepsilon}(\mathrm{Ind}_{\mathfrak a}^{\mathfrak g}S)\to V,\qquad \Pi^{\varepsilon}(x\otimes s)\mapsto (-1)^{\varepsilon |x|}xs,
    $$
    where homogeneous elements \(x\in \mathrm U(\mathfrak g)\) and \(s\in S\).

    We claim that \(\pi\) is an isomorphism.
    Since $\pi\bigl(\Pi^\varepsilon(1\otimes u)\bigr)=u$ and $u$ generates $V$, the map \(\pi\) is surjective.
    If \(\ker\pi\neq0\), Lemma \ref{nn}, applied after the parity shift, gives $\ker\pi\cap\Pi^\varepsilon S\neq0$. However, the restriction of \(\pi\) to \(\Pi^\varepsilon S\) is injective. Hence, $\ker\pi\cap\Pi^\varepsilon S=0$, a contradiction.
    Thus, $\pi$ is an isomorphism. 
    Since \(V\) is simple, Lemma \ref{simeq} implies that \(S\) is a simple \(\mathfrak s\)-module. Finally, Lemma \ref{eeqq} gives $\mathcal N(V)\cong\Pi^\varepsilon S$ as $\mathfrak s$-supermodules, where $S$ and $\mathfrak s$ are regarded as purely even supermodule and Lie superalgebra, respectively. Since $S$ is purely even, the parity of $\mathcal N(V)$ uniquely determines $\varepsilon$. Thus, $S$ and $\varepsilon$ are uniquely determined up to isomorphism.
\end{proof}
In the rest of this subsection, we shall use Ore extensions associated with $\mathfrak s$ to classify the simple $\mathfrak s$-modules $S$ defined in Theorem \ref{cc}. 

Set
\[
X=\bar f,\quad Y=-\frac12h,
\quad \Lambda=\mathrm U(\mathfrak s).
\]
The defining relation of \(\mathfrak s\) can be written as $YX-XY=X$. Thus,
\[
\Lambda\cong\mathbb C[X]\left[Y;X\frac{d}{dX}\right].
\]
We denote the first Weyl algebra
\[
\mathcal{A}_1
=
\mathbb C\langle X,\partial
\mid \partial X-X\partial=1\rangle,
\]
generated by \(X\) and \(\partial\). The assignment $Y\mapsto X\partial$ defines an embedding of $\Lambda$ into $\mathcal A_1$.
Set
$$
H=\partial X,\qquad D=\mathbb C[H],\qquad D_0=D\setminus\{0\},
$$
and define
$$
B=D_0^{-1}\mathcal A_1\cong \mathbb C(H)[X,X^{-1};\sigma],
\quad \sigma(H)=H-1,
$$
where $\sigma$ is the automorphism of $\mathbb C(H)$ determined by $\sigma(H)=H-1$. 
Recall that \cite{Bavula2020}, for nonzero polynomials \(\alpha,\beta\in D\), write \(\alpha\prec\beta\) if there are no roots \(\lambda\) and \(\mu\) of the polynomial $\alpha$ and \(\beta\) respectively and such that $\lambda-\mu$ is non-negative integral. An element
\[
\beta=\partial^m\beta_{-m}(H)+\cdots+\beta_0(H)\in\mathcal A_1,
\]
where \(m>0\) and \(\beta_i\in D\), is called \(l\)-normal if $ \beta_0\prec\beta_{-m}$ and $\beta_0\prec H$.
We then obtain the following classification result.
\begin{theorem}\label{Bavula}
    Suppose that $S$ is a simple $\mathfrak s$-module on which $X$ acts injectively. Then there exists an element
    $$
    \beta=\partial^m \beta_{-m}(H)+\dots +\beta_0(H)\in \mathcal{A}_1,
    $$
    with $m>0$ and $\beta _i(H)\in D$, that is $l$-normal and irreducible in $B$ and 
    $$
    \Lambda=\Lambda X+(\Lambda \cap B\beta ),\qquad S\cong S_{\beta}:=\Lambda/(\Lambda\cap B\beta).
    $$
    Conversely, every \(l\)-normal element \(\beta\in \mathcal A_1\) of the above form that is irreducible in \(B\) and satisfies \(\Lambda=\Lambda X+(\Lambda\cap B\beta)\) gives rise to a simple \(\Lambda\)-module \(S_\beta\) on which \(X\) acts injectively.
    
    Moreover, for modules \(S_\alpha\) and \(S_\gamma\) obtained in this way, $S_\alpha\cong S_\gamma$ if and only if $B/B\alpha\cong B/B\gamma$ as left \(B\)-modules.
\end{theorem}
\begin{proof}
    Since $\Lambda\cong \mathbb C[X]\left[Y;X\frac{d}{dX}\right]$, it is an Ore extension with \(f(X)=X\), as defined in \cite{Bavula2020}. Since $X$ is injective, the $\mathfrak s$-submodule $XS\neq 0$. By simplicity of $S$, we have $XS=S$. Hence, \(X\) acts bijectively on \(S\), so \(S\) is an \(X\)-torsion free simple $\Lambda$-module and satisfies the hypotheses in Bavula's theorem \cite[Theorem 2.10]{Bavula2020}.
    
    Therefore, there exists an element \(\beta \in \mathcal{A}_1\) of the prescribed form, which is \(l\)-normal and irreducible in \(B\), such that $\Lambda=\Lambda X+(\Lambda\cap B\beta )$ and $S\cong S_{\beta }\cong\Lambda/(\Lambda\cap B\beta ).$ Conversely, the same theorem ensures that every such element \(\beta\) determines an \(X\)-torsionfree simple \(\Lambda\)-module \(S_\beta\). The isomorphism criterion also follows from \cite[Theorem 2.10]{Bavula2020}.
\end{proof}
\begin{theorem}\label{3.38}
    Let \(V\) be a simple Whittaker \(\mathfrak g\)-supermodule of type \(\zeta=0\), and suppose that \(z_1\) and $z_2$ act on \(V\) by the scalars $\xi_1$ and $\xi_2$, respectively. Then the following statements hold:
    \begin{itemize}
        \item[(1).] If $\xi_1\neq 0$, then there exist \(\lambda,\mu\in\mathbb C\) and $\varepsilon\in \{0,1\}$, with \(\mu^2=\xi_1,\lambda=\xi_2\mu^{-1}-1\), such that \(V\cong \Pi^{\varepsilon}\mathcal L(\lambda,\mu)\).
        \item[(2).] If $\xi_1=0$ and $\xi_2\neq 0$, then no simple Whittaker $\mathfrak g$-supermodule of type $\zeta=0$ exists.
        \item[(3).] If $\xi_1=\xi_2=0$ and \(\mathfrak pV=0\), then there exist $\lambda\in \mathbb C$ and $\varepsilon\in \{0,1\}$ such that $V\cong \Pi^{\varepsilon}\mathcal L(\lambda)$.
        \item[(4).] If $\xi_1=\xi_2=0$ and \(\mathfrak pV\neq0\), then there exists \(\varepsilon\in\{0,1\}\) and an element
        \[
            \beta=\partial^m\beta_{-m}(H)+\cdots+\beta_0(H)\in \mathcal{A}_1,
        \]
        where \(m>0\) and \(\beta_i(H)\in D\), which is $l$-normal and irreducible in $B$, such that $\Lambda=\Lambda X+(\Lambda\cap B\beta)$ and $V\cong\Pi^\varepsilon \mathrm{Ind}_{\mathfrak a}^{\mathfrak g}S_{\beta}$ as $\mathfrak g$-supermodules, where $S_{\beta}=\Lambda/(\Lambda\cap B\beta )$. Conversely, for every element \(\beta\) satisfying the above conditions and every \(\varepsilon\in\{0,1\}\), the induced supermodule \(\Pi^\varepsilon \mathrm{Ind}_{\mathfrak a}^{\mathfrak g}S_{\beta}\) is a simple Whittaker \(\mathfrak g\)-supermodule of type \(\zeta=0\), on which \(z_1\) and \(z_2\) act by zero and \(\mathfrak p\) acts nontrivially. The $\mathfrak s$-module \(S_\beta\) is determined up to isomorphism by \(\mathcal N(V)\) after forgetting its parity. In particular, 
        $$
        \Pi^\varepsilon \mathrm{Ind}_{\mathfrak a}^{\mathfrak g}S_{\beta}\cong \Pi^\delta \mathrm{Ind}_{\mathfrak a}^{\mathfrak g}S_{\gamma}\iff \varepsilon=\delta,\quad B/B\beta\cong B/B\gamma.
        $$
    \end{itemize}
\end{theorem}
\begin{proof}
    The assertion follows from Remark \ref{428}, Lemmas \ref{eeqq} and \ref{simeq}, Theorems \ref{high}, \ref{cc} and \ref{Bavula}.
\end{proof}
\bigskip

\textbf{Declaration on the Use of Generative AI:} The authors used OpenAI’s ChatGPT
to improve the language and presentation of the manuscript and as an additional tool to verify the logical consistency of some proofs. All mathematical arguments, statements, and references were independently verified by the authors, who take full responsibility for the content of the paper.
\bigskip

\textbf{Data Availability Statement:} Not applicable
\bigskip

\textbf{Conflicts of Interest:} The authors have no conflicts of interest to declare that are relevant to this article.

\end{document}